\documentclass[reqno]{amsart}

\usepackage{hyperref}
\usepackage{amssymb,dsfont}

\newcommand{\intn}[2]{\ensuremath{\{  #1  ,\ldots,  #2 \}}}
\newcommand{\gene}[1]{\ensuremath{\langle #1 \rangle}}

\def\R{\mathbb{R}}
\def\N{\mathbb{N}}

\def\hau{\mathcal{H}}
\def\leb{\mathcal{L}}
\newcommand{\interieur}[1]{\ensuremath{\operatorname{int} #1}}

\def\trib{\mathcal{F}}
\def\tribbis{\mathcal{G}}
\def\prob{\mathbb{P}}
\def\esp{\mathbb{E}}
\def\coup{\mathcal{C}}

\def\dd{\mathrm{d}}

\newcommand{\diam}[1]{\ensuremath{| #1 |}}

\def\eps{\varepsilon}
\def\ph{\varphi}

\def\tree{\mathcal{T}}
\def\etiq{\mathcal{U}}
\def\raci{\varnothing}
\def\ind{\mathds{1}}

\def\calZ{\mathcal{Z}}
\def\calX{\mathcal{X}}
\def\calW{\mathcal{W}}
\def\Mu{\mathrm{M}}

\def\as{\mbox{a.s.}}

\renewcommand{\theenumi}{\alph{enumi}}

\theoremstyle{plain}
\newtheorem{thm}{Theorem}
\newtheorem{prp}{Proposition}
\newtheorem{lem}[prp]{Lemma}

\theoremstyle{definition}

\theoremstyle{remark}
\newtheorem{rem}{Remark}

\title{Random fractals and tree-indexed Markov chains}
\author{Arnaud Durand}
\address{Laboratoire d'Analyse et de Math\'ematiques Appliqu\'ees, Universit\'e Paris XII, 61 av. du G\'en\'eral de Gaulle, 94010 Cr\'eteil Cedex, France.}
\email{a.durand@univ-paris12.fr}

\begin{document}

\begin{abstract}
We study the size properties of a general model of fractal sets that are based on a tree-indexed family of random compacts and a tree-indexed Markov chain. These fractals may be regarded as a generalization of those resulting from the Moran-like deterministic or random recursive constructions considered by various authors. Among other applications, we consider various extensions of Mandelbrot's fractal percolation process.
\end{abstract}

\maketitle

\section{Introduction}\label{arbresintro}

The purpose of this paper is to study the size properties of random fractal sets based on a tree-indexed family of random compacts and a tree-indexed Markov chain. In some sense, such sets generalize the fractal sets resulting from the {\em random recursive constructions} examined by S. Graf \cite{Graf:1987fk}, R.D. Mauldin and S. Williams \cite{Mauldin:1986ei} and also K. Falconer \cite{Falconer:1986le}, which are themselves the randomized version of the recursive constructions first introduced by P. Moran \cite{Moran:1946fj} and then systematically studied by J. Hutchinson \cite{Hutchinson:1981kx}. Before presenting the fractal sets that we consider throughout the paper, we recall the main results concerning those associated with a recursive construction.

\medskip

Actually, a recursive construction is a family of compact sets indexed by the $m$-ary tree for some integer $m\geq 2$. Formally, the $m$-ary tree is the set
\[
\tree_{m}=\{\raci\}\cup\bigcup_{j=1}^\infty \intn{1}{m}^j
\]
formed by the empty word $\raci$ and the words $u=u_{1}\ldots u_{j}$ of length $j\geq 1$ in the alphabet $\intn{1}{m}$. The length $j$ of such a word $u$ is denoted by $\gene{u}$ and is called the generation of $u$. By convention, $\gene{\raci}=0$. Moreover, for any $u$ in $\tree_{m}^*=\tree_{m}\setminus\{\raci\}$, the word $\pi(u)=u_{1}\ldots u_{\gene{u}-1}$ is called the father of $u$. So, $u$ has exactly $m$ sons, which are the words $u1,\ldots,um$. Following the terminology of graph theory, the directed graph with vertex set $\tree_{m}$ and with arcs $(\pi(u),u)$, for $u\in\tree_{m}^*$, is a tree rooted at $\raci$.

Let us consider a nonempty compact subset $J_{\raci}$ of $\R^d$ (with $d\geq 1$) equal to the closure of its interior and, for any vertex $u\in\tree_{m}^*$, let us consider a compact set $J_{u}$ that is geometrically similar to $J_{\raci}$. Then, the family $(J_{u})_{u\in\tree_{m}}$ is called a recursive construction if for any vertex $u\in\tree_{m}$, the compacts indexed by the sons of $u$ are included in the compact indexed by $u$ and have disjoint interiors. The fractal set associated with such a recursive construction is
\[
K=\bigcap_{j=0}^\infty\downarrow\bigcup_{u\in\tree_{m} \atop \gene{u}=j} J_{u}.
\]
As we shall detail below, the size properties of this compact set $K$ depend mainly on the contraction ratios
\begin{equation}\label{defLu}
L_{u}=\frac{\diam{J_{u}}}{\diam{J_{\pi(u)}}},
\end{equation}
for $u\in\tree_{m}^*$, where $\diam{\cdot}$ denotes diameter.

Under the assumption that all the vectors $(L_{u1},\ldots,L_{um})$, for $u\in\tree_{m}$, are the same, P. Moran \cite{Moran:1946fj} and J. Hutchinson \cite{Hutchinson:1981kx} established that the Hausdorff dimension $s$ of $K$ satisfies
\begin{equation}\label{moraneq}
{L_{1}}^s+\ldots+{L_{m}}^s=1.
\end{equation}

Later on, K. Falconer \cite{Falconer:1986le}, S. Graf \cite{Graf:1987fk}, R.D. Mauldin and S. Williams \cite{Mauldin:1986ei} considered the case in which the compacts $J_{u}$ forming the recursive construction are random and some of them may be empty. Under the main assumption that the vectors $(L_{u1},\ldots,L_{um})$, for $u\in\tree_{m}$ such that $J_{u}$ is nonempty, are independent and identically distributed, they studied the probability that the compact set $K$ is nonempty and proved that, conditional on the fact that $K$ is nonempty, its Hausdorff dimension is almost surely equal to the infimum of all $s\geq 0$ such that
\[
\esp\left[ {L_{1}}^s+\ldots+{L_{m}}^s \right]\leq 1.
\]

A well-known example of random recursive construction is supplied by the fractal percolation process introduced by B. Mandelbrot, see \cite{Mandelbrot:1983lr}. It is defined as follows. Let us consider an integer $c\geq 2$ and a real number $p\in (0,1)$. To begin with, the square $G_{0}=[0,1]^2$ is colored black. Moreover, it may be subdivided into $c^2$ adjacent closed squares with edge length $1/c$. Each of these squares is independently colored black, with probability $p$, or white, with probability $1-p$, and the black squares form a compact set $G_{1}$. This procedure is repeated on the squares composing $G_{1}$ and the black subsquares thus obtained form a compact subset $G_{2}$ of $G_{1}$. The program is then iterated ad infinitum and yields a nested sequence $(G_{j})_{j\geq 0}$ of compact sets composed of black squares. The results of K. Falconer, S. Graf, R.D. Mauldin and S. Williams then enable to establish that the intersection over all $j\geq 0$ of the compacts $G_{j}$ is nonempty with positive probability if and only if $p>1/c^2$ and that, conditional on the fact that it is nonempty, the Hausdorff dimension of this intersection is almost surely equal to $2+\log p/\log c$ (which is the infimum of all $s\geq 0$ such that $\esp[N\,c^{-s}]\leq 1$, where $N$ is a binomial random variable with parameters $c^2$ and $p$). This result was also obtained by J. Chayes, L. Chayes and R. Durrett \cite{Chayes:1988fk} and is exposed in \cite{Chayes:1995qy,Falconer:2003oj,Graf:1988lr} too. In Section \ref{arbresapplic} below, among other applications, we explain how our results enable to study the size properties of various generalizations of Mandelbrot's fractal percolation process.

We refer to \cite{Falconer:1986le,Graf:1988lr,Lyons:2005fc,Mauldin:1986ei} for other examples of random recursive constructions. A noteworthy one is related to the zero set of the Brownian bridge. Indeed, the results of K. Falconer, S. Graf, R.D. Mauldin and S. Williams may be applied in order to recover a result of S. Taylor \cite{Taylor:1955qy} according to which the Hausdorff dimension of this set is almost surely equal to $1/2$, see \cite{Graf:1988lr} for details.

\medskip

From now on and except in Section \ref{arbresapplic} (in which we give several applications of our results), we shall always consider the general case.

Various refinements and extensions of the aforementioned results were obtained. To begin with, S. Graf, R.D. Mauldin and S. Williams \cite{Graf:1988lr} found a gauge function $h$ satisfying $0<\hau^h(K)<\infty$ with probability one, conditional on the fact that $K$ is nonempty, where $\hau^h$ denotes the Hausdorff $h$-measure (see \cite{Rogers:1970wb} for the definition). Later, A. Berlinkov and R.D. Mauldin \cite{Berlinkov:2003kx,Berlinkov:2002fj} studied the packing dimension and measures of the compact $K$ (see \cite{Falconer:2003oj} for the definitions). Moreover, several authors established the existence of self-similar random measures carried by compacts analogous to $K$ and performed their multifractal analysis, see \cite{Arbeiter:1990lr,Falconer:1994fk,Hutchinson:1998uq,Olsen:1994qy}. Let us also mention that Y. Pesin and H. Weiss \cite{Pesin:1996bb}, and also Y. Kifer \cite{Kifer:1995nu,Kifer:1996yq}, employed some techniques from the theory of dynamical systems with a view to determining the Hausdorff dimension of various random sets that are built recursively.

In all the works cited above, the vectors $(L_{u1},\ldots,L_{um})$ giving the contraction ratios across generations are independent and identically distributed. A. Dryakhlov and A. Tempelman \cite{Dryakhlov:2001vn} proposed a way to relax this assumption. Specifically, under the main assumption that the vector $(L_{u1},\ldots,L_{um})$ associated with a given vertex $u\in\tree_{m}$ is correlated with the vectors corresponding to a fixed number of ancestors of $u$, they established that the Hausdorff dimension of the compact $K$ is almost surely equal to a specific value which can be computed in terms of the distributions of the contraction ratios. Note that this phenomenon could be anticipated thanks to Kolmogorov's zero-one law.

\medskip

In the same vein, Y.-Y. Liu, Z.-Y. Wen and J. Wu \cite{Liu:2004rt} introduced another generalization of the previous random recursive constructions, in which the vectors giving the contraction ratios need not be identically distributed. This generalization is the first step towards the fractal sets that we study in this paper, so we now recall it. We begin by replacing the $m$-ary tree $\tree_{m}$ by the tree $\etiq_{0}$ in which every vertex with generation $j$ has exactly $m_{j}$ sons, where $(m_{j})_{j\geq 0}$ is a given sequence of integers greater than one. More precisely,
\begin{equation}\label{defetiq0}
\etiq_{0}=\left\{ u\in\etiq \:\bigl|\: \forall j\in\intn{1}{\gene{u}} \quad u_{j}\leq m_{j-1} \right\},
\end{equation}
where $\etiq$ denotes the set formed by the empty word $\raci$ and the words of finite length in the alphabet $\N=\{1,2,\ldots\}$. Of course, if $m_{j}=m$ for every integer $j\geq 0$ and some integer $m\geq 2$, then $\etiq_{0}$ is just the $m$-ary tree $\tree_{m}$.

We then consider a random family $(J_{u})_{u\in\etiq_{0}}$ of compact subsets of $\R^d$ indexed by the tree $\etiq_{0}$, having positive diameter and satisfying the following properties:
\renewcommand{\theenumi}{\Alph{enumi}}
\begin{enumerate}
\item\label{condcompact1} For any vertex $u\in\etiq_{0}$, the compacts $J_{uk}$, for $k\in\intn{1}{m_{\gene{u}}}$, are subsets of $J_{u}$ with disjoint interiors.
\item\label{condcompact2} There exists a real $\kappa>0$ such that for any $u\in\etiq_{0}$, the Lebesgue measure of the interior of $J_{u}$ is at least $\kappa\diam{J_{u}}^d$.
\end{enumerate}
We also make an assumption about the distribution of the contraction ratios defined by (\ref{defLu}). Specifically, let $\underline{\beta}$ and $\overline{\beta}$ be two real numbers enjoying $0<\underline{\beta}\leq\overline{\beta}<1$ and, for any $j\geq 0$, let $\mu_{j}$ be a probability measure on $[\underline{\beta},\overline{\beta}]^{m_{j}}$. We suppose that:
\begin{enumerate}\setcounter{enumi}{2}
\item\label{condcompact3} The $(L_{uk})_{k\in\intn{1}{m_{\gene{u}}}}$, for $u\in\etiq_{0}$, are independent random vectors with distribution $\mu_{\gene{u}}$.
\end{enumerate}
Under assumptions somewhat stronger than (\ref{condcompact1}-\ref{condcompact3}), Y.-Y. Liu, Z.-Y. Wen and J. Wu studied the size properties of the random compact set
\[
K=\bigcap_{j=0}^\infty\downarrow\bigcup_{u\in\etiq_{0}\atop\gene{u}=j} J_{u}.
\]
To be specific, they established that its Hausdorff dimension is almost surely equal to a certain value which can be expressed in terms of the probability measures $\mu_{j}$. Note that this compact set may also be obtained in the following manner. Assertion (\ref{condcompact3}) ensures that with probability one, for any sequence $\zeta=(\zeta_{j})_{j\geq 1}$ in $\N$ enjoying $\zeta_{1}\ldots\zeta_{j}\in\etiq_{0}$ for all $j\geq 1$, the diameter $\diam{J_{\zeta_{1}\ldots\zeta_{j}}}$ is at most $\diam{J_{\raci}}\overline{\beta}^j$, so that it tends to zero as $j\to\infty$. Hence, there exists a unique point $x_{\zeta}$ in $\R^d$ such that
\begin{equation}\label{defxzeta}
\{x_{\zeta}\}=\bigcap_{j=1}^\infty\downarrow J_{\zeta_{1}\ldots\zeta_{j}}.
\end{equation}
A standard diagonal argument then shows that the compact $K$ is also equal to the collection of all such points $x_{\zeta}$.

\medskip

The fractal set that we study below is a random subset $\Theta$ of $K$ chosen according to a tree-indexed Markov chain which we now introduce. Let us consider a family $(X_{u})_{u\in\etiq_{0}}$ of $\{0,1\}$-valued random variables which is independent of the family $(L_{u})_{u\in\etiq_{0}^*}$, where $\etiq_{0}^*=\etiq_{0}\setminus\{\raci\}$ (a case in which these two families need not be independent is briefly discussed in Section \ref{arbresgeneralisation}). In addition, for any $j_{0}\geq 1$, let $\etiq_{j_{0}}$ denote the tree obtained by replacing the sequence $(m_{j})_{j\geq 0}$ by the sequence $(m_{j_{0}+j})_{j\geq 0}$ in the definition (\ref{defetiq0}) of $\etiq_{0}$ and let $\etiq_{j_{0}}^*=\etiq_{j_{0}}\setminus\{\raci\}$. Note that, for any vertex $u\in\etiq_{0}$, the set $u\etiq_{\gene{u}}$ (i.e.~the set of all concatenations of the word $u$ with words of $\etiq_{\gene{u}}$) is the subtree of $\etiq_{0}$ which is rooted at $u$. Hence, the $\sigma$-field
\[
\tribbis_{u}=\sigma( X_{v},\ v\in\etiq_{0}\setminus (u\etiq_{\gene{u}}^*) )
\]
can be seen as the past before $u$ in the tree $\etiq_{0}$. Indeed, this $\sigma$-field is generated by the random variables corresponding to the vertices of $\etiq_{0}$ which are not descended from $u$. Conversely, the future after $u$ begins with its sons, which are the vertices $uk$ for $k\in\intn{1}{m_{\gene{u}}}$.

For every integer $j\geq 0$ and every $t\in\{0,1\}$, let $\nu_{t,j}$ be a probability measure on $\{0,1\}^{m_{j}}$. From now on and except in Section \ref{arbresgeneralisation}, we assume that the process $(X_{u})_{u\in\etiq_{0}}$ is a Markov chain with transition probability measures $\nu_{t,j}$ for $j\geq 0$ and $t\in\{0,1\}$, which means that the following {\em Markov condition} holds:
\begin{enumerate}\setcounter{enumi}{3}
\item\label{condMarkov} For any vertex $u\in\etiq_{0}$ and any subset $A$ of $\{0,1\}^{m_{\gene{u}}}$,
\[
\prob\bigl( (X_{uk})_{k\in\intn{1}{m_{\gene{u}}}} \in A \:\bigl|\: \tribbis_{u} \bigr)=\nu_{X_{u},\gene{u}}(A).
\]
\end{enumerate}
\renewcommand{\theenumi}{\alph{enumi}}
Informally, for every vertex $u\in\etiq_{0}$, the vector $(X_{uk})_{k\in\intn{1}{m_{\gene{u}}}}$ depends only on the value of $X_{u}$ and the generation $\gene{u}$, conditionally on the past before $u$. The definition of a tree-indexed Markov chain that we adopt here may be compared with that introduced by I. Benjamini and Y. Peres \cite{Benjamini:1994yv}. Actually, the Markov chains that they considered correspond to the particular case in which the measures $\nu_{t,j}$, for $t\in\{0,1\}$, are the products ${\lambda_{t}}^{\otimes m_{j}}$ for some fixed probability measure $\lambda_{t}$ on $\{0,1\}$. However, note that we restrict our attention to $\{0,1\}$-valued Markov chains and trees of the form $\etiq_{0}$, whereas I. Benjamini and Y. Peres did not.

An important consequence of the Markov condition (\ref{condMarkov}) is that for any integer $j_{0}\geq 0$, conditionally on the $\sigma$-field generated by the variables $X_{w}$ for $w\in\etiq_{0}$ with $\gene{w}\leq j_{0}$, the processes $(X_{uv})_{v\in\etiq_{j_{0}}}$, for $u\in\etiq_{0}$ with generation $j_{0}$, are independent Markov chains with transition probability measures $\nu_{t,j_{0}+j}$, for $j\geq 0$ and $t\in\{0,1\}$.

To obtain the set $\Theta$, we only keep in the compact set $K$ the points $x_{\zeta}$ resulting via (\ref{defxzeta}) from a sequence $\zeta=(\zeta_{j})_{j\geq 1}$ in $\N$ enjoying $X_{\zeta_{1}\ldots\zeta_{j}}=1$ for all $j$ large enough. To be specific, for any vertex $u\in\etiq_{0}$, let
\begin{equation}\label{deftauu}
\tau_{u}=\left\{ v\in u\etiq_{\gene{u}} \:\bigl|\: \forall j\in\intn{\gene{u}}{\gene{v}} \quad X_{v_{1}\ldots v_{j}}=1 \right\}.
\end{equation}
If $X_{u}=0$, then the set $\tau_{u}$ is empty. Otherwise, $\tau_{u}$ is the largest subtree of $\etiq_{0}$ rooted at $u$ and formed by vertices mapped to the state $1$ by the Markov chain $X$. The boundary of $\tau_{u}$ is
\begin{equation}\label{deffrtauu}
\partial\tau_{u}=\left\{ \zeta=(\zeta_{j})_{j\geq 1}\in \N^{\N} \:\bigl|\: \forall j\geq\gene{u} \quad \zeta_{1}\ldots\zeta_{j}\in\tau_{u} \right\}.
\end{equation}
The points $x_{\zeta}$ resulting from all the sequences $\zeta$ of this boundary form the set
\begin{equation}\label{defKu}
K_{u}=\bigcup_{\zeta\in\partial\tau_{u}} \{x_{\zeta}\}.
\end{equation}
As observed before, a standard diagonal argument enables to show that this last set is a compact subset of $K$. The subset $\Theta$ of $K$ that we study in the following is then the $F_{\sigma}$-set
\begin{equation}\label{defTheta}
\Theta=\bigcup_{u\in\etiq_{0}} K_{u}.
\end{equation}
A point of $K$ thus also belongs to $\Theta$ if and only if it can be written on the form $x_{\zeta}$ for some sequence $\zeta=(\zeta_{j})_{j\geq 1}$ in $\N$ such that $X_{\zeta_{1}\ldots\zeta_{j}}=1$ for all $j$ large enough. Furthermore, note that the randomness in the construction of $\Theta$ lies both in the family $(J_{u})_{u\in\etiq_{0}}$ of compact sets and in the Markov chain $(X_{u})_{u\in\etiq_{0}}$.

\medskip

The fractal sets obtained by dint of the recursive constructions introduced by the aforementioned authors may actually be seen as particular cases of the set $\Theta$. Indeed, if $X_{\raci}=1$ with probability one and $\nu_{1,j}$ is the point mass at $(1,\ldots,1)$ for any $j\geq 0$, then $\tau_{\raci}$ is almost surely equal to the whole tree $\etiq_{0}$, so that $\Theta$ is almost surely equal to the whole compact set $K$. One thus recovers the fractal sets obtained through the constructions introduced by Y.-Y. Liu, Z.-Y. Wen and J. Wu. Likewise, one can obtain the fractals resulting from the constructions of K. Falconer, S. Graf, R.D. Mauldin and S. Williams by letting $\nu_{0,j}$ be the point mass at $(0,\ldots,0)$ and assuming that the probability measures $\nu_{1,j}$ and $\mu_{j}$ do not depend on $j$. Therefore, in some sense, the results of this paper can be seen as a generalization of those established in \cite{Falconer:1986le,Graf:1987fk,Liu:2004rt,Mauldin:1986ei}. They also apply when the family $(J_{u})_{u\in\etiq_{0}}$ of compacts is deterministic, as in the works of P. Moran and J. Hutchinson. In this case, the distributions $\mu_{j}$ of the contraction ratios are simply a fixed point mass.

\medskip

The rest of the paper is organized as follows. The main results are exposed in Section \ref{arbrespres} and several applications are given in Section \ref{arbresapplic}. In Section \ref{arbresprel}, we establish various ancillary lemmas: we study a particular family of branching processes in varying environment related to the underlying Markov chain, we adapt a method proposed by K. Falconer \cite{Falconer:1986le} to exhibit a connection between the size properties of the compact sets $K_{u}$ composing $\Theta$ and the question of the existence of positive flows in certain random networks and we generalize some techniques exposed by R. Lyons and Y. Peres \cite{Lyons:2005fc} concerning percolation on trees. Sections \ref{arbresmajor} and \ref{arbresminor} are respectively devoted to giving an upper and a lower bound on the Hausdorff dimension of the compact sets $K_{u}$ and the main results of the paper are proven in Section \ref{arbresmainproofs}. Lastly, Section \ref{arbresgeneralisation} briefly discusses a simple extension of these results to a case in which the families $(L_{u})_{u\in\etiq_{0}^*}$ and $(X_{u})_{u\in\etiq_{0}}$ need not be independent.

\section{Statement of the results}\label{arbrespres}

The main results of the paper concern the distribution of the Hausdorff dimension of the set $\Theta$. It is not obvious that the dimension of $\Theta$ is a random variable, i.e.~is measurable with respect to the $\sigma$-field of the underlying probability space. However, this is true because it is measurable with respect to the $\sigma$-field generated by the contraction ratios $L_{u}$, for $u\in\etiq_{0}^*$, and the states $X_{u}$, for $u\in\etiq_{0}$, which are random variables. We refer to Remark~\ref{remarkmeasurability} below for details. Note that there is no need to make measurability assumptions on the compacts $J_{u}$ themselves.

Before stating the results, let us recall the definition of Hausdorff dimension. To begin with, for every real $s>0$, the $s$-dimensional Hausdorff measure of a subset $F$ of $\R^d$ is given by
\[
\hau^s(F)=\lim_{\eps\downarrow 0}\uparrow\hau^s_{\eps}(F) \qquad\text{with}\qquad \hau^s_{\eps}(F)=\inf_{F\subseteq\bigcup_{p} U_{p}\atop\diam{U_{p}}<\eps} \sum_{p=1}^\infty \diam{U_{p}}^s,
\]
where the infimum is taken over all sequences $(U_{p})_{p\geq 1}$ of sets with $F\subseteq\bigcup_{p} U_{p}$ and $\diam{U_{p}}<\eps$ for all $p$. Note that $\hau^s$ is a Borel measure on $\R^d$, see \cite{Rogers:1970wb}. The Hausdorff dimension of a nonempty set $F\subseteq\R^d$ is then defined by
\[
\dim F=\sup\{s\in (0,d) \:|\: \hau^s(F)=\infty\}=\inf\{s\in (0,d) \:|\: \hau^s(F)=0\},
\]
with the convention that the supremum (resp.~infimum) of the empty set is zero (resp.~$d$), see \cite{Falconer:2003oj}. In addition, we agree that the empty set $\emptyset$ has Hausdorff dimension $-\infty$.

\medskip

We also need to introduce a family of real numbers $\alpha_{s,j}$ related to the the transition probability measures $\nu_{1,j}$ of the Markov chain $(X_{u})_{u\in\etiq_{0}}$ and the distributions $\mu_{j}$ of the ratios $(L_{u})_{u\in\etiq_{0}^*}$. Specifically, for any real $s$ and any integer $j\geq 0$, let
\begin{equation}\label{defalphasj}
\alpha_{s,j}=\int_{\{0,1\}^{m_{j}}}\int_{[\underline{\beta},\overline{\beta}]^{m_{j}}} \sum_{k=1}^{m_{j}} {\ell_{k}}^s x_{k} \ \mu_{j}(\dd\ell)\nu_{1,j}(\dd x),
\end{equation}
where $\ell_{1},\ldots,\ell_{m_{j}}$ and $x_{1},\ldots,x_{m_{j}}$ are the coordinates of $\ell\in [\underline{\beta},\overline{\beta}]^{m_{j}}$ and $x\in\{0,1\}^{m_{j}}$, respectively. Note that $\alpha_{s,j}$ can be seen as a generalization of the left-hand side of the Moran equation (\ref{moraneq}).

The reals $\alpha_{s,j}$ enable us to introduce a number $d_{*}\in [-\infty,\infty)$ which governs the Hausdorff dimension of $\Theta$. It is defined as follows. If
\begin{equation}\label{defunderj}
\underline{j}=\inf\{ j_{0}\geq 0 \:|\: \forall j\geq j_{0} \quad \alpha_{0,j}>0 \}
\end{equation}
is infinite, then let $d_{*}=-\infty$. Conversely, if $\underline{j}$ is finite, then let us consider the function $\rho$ defined on $\R$ by
\begin{equation}\label{defrhodim}
\rho(s)=\liminf_{j\to\infty} \frac{1}{j}\sum_{n=\underline{j}}^{j-1}\log\alpha_{s,n}.
\end{equation}
Owing to the fact that the distributions $\mu_{j}$ are supported on $[\underline{\beta},\overline{\beta}]^{m_{j}}$, the function $\rho$ is either the constant function equal to $-\infty$ on $\R$ or a decreasing bijection from $\R$ onto $\R$. In both cases, one can consider
\begin{equation}\label{defdetoile}
d_{*}=\sup\{ s\in\R \:|\: \rho(s)>0 \}=\inf\{ s\in\R \:|\: \rho(s)<0 \}.
\end{equation}
It is possible to give another expression of $d_{*}$ when $\underline{j}$ is finite. Indeed, in this case, the function $s\mapsto \log\alpha_{s,\underline{j}}+\ldots+\log\alpha_{s,j-1}$ is a bijection from $\R$ onto $\R$ for any $j>\underline{j}$, so it has a unique zero denoted by $d_{j}$. One then readily verifies that
\[
d_{*}=\liminf\limits_{j\to\infty} d_{j}.
\]

The following result, which is established in Section \ref{arbresmainproofs}, gives the possible values of the Hausdorff dimension of the set $\Theta$ defined by (\ref{defTheta}).

\begin{thm}\label{loidimTheta}
With probability one,
\[
\left\{\begin{array}{rcl}
d_{*}<0 & \qquad\Longrightarrow\qquad & \dim\Theta=-\infty \\[2mm]
d_{*}\geq 0 & \qquad\Longrightarrow\qquad & \dim\Theta\in\{-\infty,d_{*}\}.
\end{array}\right.
\]
\end{thm}

\medskip

In order to complete the description of the distribution of the Hausdorff dimension of $\Theta$, there remains to study the probability that it is equal to $-\infty$ in the case where $d_{*}$ is nonnegative. This amounts to examining the probability that the set $\Theta$ is empty. To this end, we need to introduce several notations related to the underlying Markov chain $(X_{u})_{u\in\etiq_{0}}$. For any integer $j\geq 0$, let
\[
S_{j}=\left\{ u\in\etiq_{0} \:|\: \gene{u}=j\text{ and }X_{u}=1 \right\}.
\]
This is the set of all vertices with generation $j$ that are mapped to the state $1$ by the Markov chain. Moreover, let us consider the generating function of the cardinality of $S_{j}$, that is,
\begin{equation}\label{defPhij}
\Phi_{j}:z\mapsto\esp[z^{\# S_{j}}].
\end{equation}
The generating functions $\Phi_{0},\Phi_{1},\ldots$ may easily be computed as follows in terms of the transition probability measures $\nu_{t,j}$ of the Markov chain. Actually, the Markov condition (\ref{condMarkov}) implies that for any integer $j\geq 0$ and any complex number $z$,
\begin{equation}\label{recfctgeneSj}
\esp[z^{\# S_{j+1}}\:|\: X_{u},\,\gene{u}\leq j]=\ph_{1,j}(z)^{\# S_{j}}\ph_{0,j}(z)^{m_{0}\cdot\ldots\cdot m_{j-1}-\# S_{j}},
\end{equation}
where the functions $\ph_{0,j}$ and $\ph_{1,j}$ are given by
\begin{equation}\label{defphtj}
\ph_{t,j}(z)=\int_{\{0,1\}^{m_{j}}} z^{x_{1}+\ldots+x_{m_{j}}} \nu_{t,j}(\dd x).
\end{equation}
Taking expectations in (\ref{recfctgeneSj}), it follows that the generating functions $\Phi_{0},\Phi_{1},\ldots$ may be calculated recursively using the formulas:
\begin{equation}\label{Phijrec}
\left\{\begin{array}{rcl}
\Phi_{0}(z) &=& \prob(X_{\raci}=0)+\prob(X_{\raci}=1)\cdot z \\[2mm]
\Phi_{j+1}(z) &=& \ph_{0,j}(z)^{m_{0}\cdot\ldots\cdot m_{j-1}}\cdot\Phi_{j}(\ph_{1,j}(z)/\ph_{0,j}(z)).
\end{array}\right.
\end{equation}

Alongside with that, we need to consider the sequence $(f_{j})_{j\geq 0}$ defined by
\begin{equation}\label{deffj}
\forall j\geq 0 \qquad f_{j}=\lim_{n\uparrow\infty}\uparrow \ph_{1,j}\circ\ldots\circ\ph_{1,j+n}(0).
\end{equation}
As shown by Lemma \ref{lemlienfj} below, every real $f_{j}$ is in fact the extinction probability of a branching process in varying environment related to the transition probability measures $\nu_{1,j+n}$, for $n\geq 0$.

The following result, which is proven in Section \ref{arbresmainproofs}, provides an expression of the probability that $\Theta$ is empty, thereby completing the description of the distribution of the Hausdorff dimension of this set.

\begin{thm}\label{thmprobThetavide}
If $d_{*}\geq 0$, then
\[
\prob(\Theta=\emptyset)=\lim_{j\uparrow\infty}\downarrow \Phi_{j}(f_{j}).
\]
\end{thm}

When $f_{j}$ vanishes for some $j$, it is possible to provide an expression of the probability that $\Theta$ is empty that is more explicit than that given by Theorem \ref{thmprobThetavide}. This is the purpose of the following result, which is proven in Section \ref{arbresmainproofs}.

\begin{prp}\label{altthmprobThetavide}
If $d_{*}\geq 0$ and $f_{j_{*}}=0$ for some integer $j_{*}\geq 0$, then
\[
\prob(\Theta=\emptyset)=\Phi_{j_{*}}(0)\cdot\prod_{j=j_{*}}^\infty\ph_{0,j}(0)^{m_{0}\cdot\ldots\cdot m_{j-1}}.
\]
\end{prp}

In whole generality and especially when the conditions of Proposition \ref{altthmprobThetavide} do not hold, it seems awkward to provide an expression of the probability that $\Theta$ is empty which is both tractable and more explicit than that given by Theorem \ref{thmprobThetavide}. Instead, we supply necessary and sufficient conditions for $\Theta$ to be empty with positive probability and with probability one, respectively. These conditions are expressed by means of two sequences $(\underline{\sigma}_{j})_{j\geq 0}$ and $(\overline{\sigma}_{j})_{j\geq 0}$ which are defined by
\begin{equation}\label{defsigmaj}
\left\{\begin{array}{rcl}
\underline{\sigma}_{j} &=& \max\Biggl(1,\sum\limits_{n=j}^\infty \frac{ \ph_{1,n}'(1)+\ph_{1,n}(0)-1 }{ \ph_{1,n}'(1) \prod\limits_{\ell=j}^n \ph_{1,\ell}'(1) } + \limsup\limits_{n\to\infty}\frac{1}{\prod\limits_{\ell=j}^n \ph_{1,\ell}'(1)}\Biggr) \\[6mm]
\overline{\sigma}_{j} &=& \sum\limits_{n=j}^\infty \frac{ \ph_{1,n}''(1) }{ \ph_{1,n}'(1) \prod\limits_{\ell=j}^n \ph_{1,\ell}'(1) } + \liminf\limits_{n\to\infty}\frac{1}{\prod\limits_{\ell=j}^n \ph_{1,\ell}'(1)}.
\end{array}\right.
\end{equation}
Actually, these sequences lead to a lower and an upper bound on $f_{j}$. More precisely, Lemma \ref{lemrecencfj} below shows that $1/(1-f_{j})$ is between $\underline{\sigma}_{j}$ and $\overline{\sigma}_{j}$ for every integer $j\geq 0$. Note that $f_{j}$ is thus equal to one if $\underline{\sigma}_{j}=\infty$. This is so in particular when $j$ is less than the number $\underline{j}$ defined by (\ref{defunderj}).

Let us now give a necessary and a sufficient condition for the set $\Theta$ to be empty with positive probability or with probability one, respectively. Thanks to Proposition \ref{altthmprobThetavide}, we can obviously restrict our attention to case in which $d_{*}\geq 0$ and
\begin{equation}\label{condfjpos}
\forall j\geq 0 \qquad f_{j}>0.
\end{equation}
Note that, under these assumptions, $\underline{j}$ is necessarily finite, as $d_{*}$ is chosen to be equal to $-\infty$ when $\underline{j}$ is infinite. The next two results are established in Section \ref{arbresmainproofs}.

\begin{prp}\label{condnecsufThetavide1}
Let us assume that $d_{*}\geq 0$ and that (\ref{condfjpos}) holds. Then,
\begin{eqnarray}
& & \exists j_{0}\geq 0 \qquad \sum\limits_{j=j_{0}}^\infty -m_{0}\cdot\ldots\cdot m_{j-1}\cdot\log\ph_{0,j}\left(1-\frac{1}{\underline{\sigma}_{j+1}}\right)<\infty \label{conddimTheta1} \\
\Longrightarrow & & \prob(\Theta=\emptyset)>0 \nonumber \\
\Longrightarrow & & d_{*}=0 \qquad\text{or}\qquad \sum\limits_{j=0}^\infty -m_{0}\cdot\ldots\cdot m_{j-1}\cdot\log\ph_{0,j}\left(1-\frac{1}{\overline{\sigma}_{j+1}}\right)<\infty. \label{conddimTheta1bis}
\end{eqnarray}
\end{prp}

\begin{prp}\label{condnecsufThetavide2}
Let us assume that $d_{*}\geq 0$ and that (\ref{condfjpos}) holds. Then,
\begin{eqnarray}
& & \underline{\sigma}_{\underline{j}}=\infty \qquad\text{or}\qquad [\Phi_{\underline{j}}(0)=1 \quad\text{and}\quad \forall j\geq\underline{j} \quad \ph_{0,j}(0)=1] \label{conddimTheta2} \\
\Longrightarrow & & \prob(\Theta=\emptyset)=1 \nonumber \\
\Longrightarrow & & \overline{\sigma}_{\underline{j}}=\infty \qquad\text{or}\qquad [\Phi_{\underline{j}}(0)=1 \quad\text{and}\quad \forall j\geq\underline{j} \quad \ph_{0,j}(0)=1]. \label{conddimTheta3}
\end{eqnarray}
\end{prp}

\begin{rem}\label{remarksigma}
For particular choices of the transition probability measures $\nu_{1,j}$ of the Markov chain, the sequences $(\underline{\sigma}_{j})_{j\geq 0}$ and $(\overline{\sigma}_{j})_{j\geq 0}$ behave comparably, so that Propositions \ref{condnecsufThetavide1} and \ref{condnecsufThetavide2} actually provide a criterion to know if $\Theta$ is empty with positive probability or with probability one respectively, see Proposition \ref{corremarksigma} below for an illustration of this remark.
\end{rem}

\begin{rem}
An important consequence of Theorem \ref{loidimTheta} and Proposition \ref{condnecsufThetavide1} is the fact that, under the assumptions that $d_{*}>0$ and that (\ref{condfjpos}) holds,
\[
\sum_{j=0}^\infty -m_{0}\cdot\ldots\cdot m_{j-1}\cdot\log\ph_{0,j}\left(1-\frac{1}{\overline{\sigma}_{j+1}}\right)=\infty \quad\Longrightarrow\quad \as \quad \dim\Theta=d_{*}.
\]
\end{rem}

\section{Applications}\label{arbresapplic}

\subsection{Generalizations of Mandelbrot's fractal percolation process}

Various generalizations of the fractal percolation process which we briefly described in Section \ref{arbresintro} were considered, see \cite{Chayes:1995qy} and the references therein. In what follows, we introduce several new ones, to which the results of Section \ref{arbrespres} may be applied.

\subsubsection{A generalization of the Bernoulli case}

Let $(c_{j})_{j\geq 0}$ denote a bounded sequence of integers greater than one, let $\pi\in [0,1]$ and let $(p_{j})_{j\geq 0}$ and $(q_{j})_{j\geq 0}$ be two sequences in $[0,1]$. To begin with, the cube $[0,1]^d$ is colored black, with probability $\pi$, or white, with probability $1-\pi$. Moreover, it is the union of ${c_{0}}^d$ adjacent closed subcubes with edge length $1/c_{0}$. If $[0,1]^d$ is black (resp.~white), then each of these subcubes is independently colored black, with probability $p_{0}$ (resp.~$q_{0}$), or white, with probability $1-p_{0}$ (resp.~$1-q_{0}$). Each of the black (resp.~white) subcubes is itself the union of ${c_{1}}^d$ adjacent closed cubes with edge length $1/(c_{0}c_{1})$ and each of these last cubes is then independently colored black, with probability $p_{1}$ (resp.~$q_{1}$), or white, with probability $1-p_{1}$ (resp.~$1-q_{1}$). This program is iterated ad infinitum. Note that the usual percolation process corresponds to the case in which $\pi=1$, the sequence $(c_{j})_{j\geq 0}$ is constant, the sequence $(p_{j})_{j\geq 0}$ is constant and equal to a given $p\in (0,1)$ and the sequence $(q_{j})_{j\geq 0}$ is equal to zero. We are interested in the distribution of the Hausdorff dimension of the set of points belonging to black cubes from a certain stage onward.

More formally, this set corresponds to the set $\Theta$ given by (\ref{defTheta}) when the tree $\etiq_{0}$ is defined by (\ref{defetiq0}) for $m_{j}={c_{j}}^d$ and the family of compact sets $(J_{u})_{u\in\etiq_{0}}$ and the Markov chain $(X_{u})_{u\in\etiq_{0}}$ are as follows. For any integer $j\geq 0$, let us consider a bijection $k\mapsto (r_{j,1}(k),\ldots,r_{j,d}(k))$ from $\intn{1}{m_{j}}$ onto $\intn{0}{c_{j}-1}^d$. Then, let $J_{\raci}=[0,1]^d$ and, for any vertex $u\in\etiq_{0}^*$, let
\begin{equation}\label{defJuperco}
J_{u}=\left(\sum_{j=1}^{\gene{u}}\frac{r_{j,1}(u_{j})}{c_{0}\cdot\ldots\cdot c_{j-1}},\ldots,\sum_{j=1}^{\gene{u}}\frac{r_{j,d}(u_{j})}{c_{0}\cdot\ldots\cdot c_{j-1}}\right)+\frac{1}{c_{0}\cdot\ldots\cdot c_{\gene{u}-1}} [0,1]^d.
\end{equation}
Observe that Assertions (\ref{condcompact1})-(\ref{condcompact3}) hold when the probability measure $\mu_{j}$ is the point mass at $(1/c_{j},\ldots,1/c_{j})$ for every integer $j\geq 0$. Moreover, let $(X_{u})_{u\in\etiq_{0}}$ be a $\{0,1\}$-valued Markov chain such that $X_{\raci}=1$ with probability $\pi$ and with transitions given by the product measures
\begin{equation}\label{kernelbernoulli}
\left\{\begin{array}{rcl}
\nu_{0,j} & = & (q_{j}\delta_{1}+(1-q_{j})\delta_{0})^{\otimes m_{j}} \\[2mm]
\nu_{1,j} & = & (p_{j}\delta_{1}+(1-p_{j})\delta_{0})^{\otimes m_{j}},
\end{array}\right.
\end{equation}
where $\delta_{x}$ denotes the point mass at $x$. Then, the black (resp.~white) cubes correspond to the compacts $J_{u}$ indexed by vertices $u\in\etiq_{0}$ enjoying $X_{u}=1$ (resp.~$0$). Thus, as announced previously, the set $\Theta$ is the set of points that belong to black cubes from some stage onward.

With a view to applying the results of Section \ref{arbrespres}, we need to consider the infimum $\underline{j}$ of all integers $j_{0}\geq 0$ such that $p_{j}$ is positive for every $j\geq j_{0}$ and then to let $d_{*}=-\infty$ if $\underline{j}$ is infinite and
\[
d_{*}=d+\liminf_{j\to\infty}\frac{\log p_{\underline{j}}+\ldots+\log p_{j}}{\log c_{\underline{j}}+\ldots+\log c_{j}}
\]
otherwise. If $d_{*}$ is negative, then the set $\Theta$ is almost surely empty by virtue of Theorem \ref{loidimTheta}. Conversely, if $d_{*}$ is nonnegative, then Theorem \ref{loidimTheta} ensures that with probability one, the set $\Theta$ either is empty or has Hausdorff dimension $d_{*}$. Moreover, Propositions \ref{altthmprobThetavide}, \ref{condnecsufThetavide1} and \ref{condnecsufThetavide2} provide additional information about the probability that $\Theta$ is empty in the case where $d_{*}$ is nonnegative.

Remark that if $p_{j}$ tends to one as $j\to\infty$, then $d_{*}$ is equal to $d$, so that the set $\Theta$ has Hausdorff dimension $d$ with probability one when it is nonempty. Let us consider the particular case in which $d$ is equal to two, $\pi$ is equal to one, the sequence $(c_{j})_{j\geq 0}$ is constant and the sequence $(q_{j})_{j\geq 0}$ is equal to zero. Then, the set $\Theta$ has Hausdorff dimension two with probability one when it is nonempty. This observation motivated the authors of \cite{Chayes:1997lr} to use this set as a better mean field approximation of the planar Brownian path, which has Hausdorff dimension two, than the sets resulting from the usual fractal percolation process, which necessarily have Hausdorff dimension less than two.

\subsubsection{A generalization of the microcanonical case}

The difference with the previous case is that the transition probability measures defined by (\ref{kernelbernoulli}) are replaced by
\[
\nu_{0,j}=\frac{1}{\binom{m_{j}}{b_{j}}}\sum_{x\in\{0,1\}^{m_{j}}\atop \sum_{k} x_{k}=b_{j}}\delta_{x}
\qquad\text{and}\qquad
\nu_{1,j}=\frac{1}{\binom{m_{j}}{a_{j}}}\sum_{x\in\{0,1\}^{m_{j}}\atop \sum_{k} x_{k}=a_{j}}\delta_{x},
\]
where $a_{j},b_{j}\in\intn{0}{m_{j}}$ for any integer $j\geq 0$. Informally, the cubes are now colored as follows. As in the previous situation, the cube $[0,1]^d$ is colored black, with probability $\pi$, or white, with probability $1-\pi$. If $[0,1]^d$ is black (resp.~white), then among its $m_{0}={c_{0}}^d$ adjacent closed subcubes with edge length $1/c_{0}$, exactly $a_{0}$ (resp.~$b_{0}$) are colored black and their positions are chosen uniformly. Next, if a cube with edge length $1/c_{0}$ is black (resp.~white), then exactly $a_{1}$ (resp.~$b_{1}$) of its $m_{1}={c_{1}}^d$ adjacent closed subcubes with edge length $1/(c_{0}c_{1})$ are colored black and their positions are chosen uniformly again. The above program is then reenacted ad infinitum. Note that the usual microcanonical fractal percolation process corresponds to the case in which $\pi=1$, the sequences $(c_{j})_{j\geq 0}$ and $(a_{j})_{j\geq 0}$ are constant and the sequence $(b_{j})_{j\geq 0}$ is equal to zero.

The results of Section \ref{arbrespres} prompt us to consider the infimum $\underline{j}$ of all integers $j_{0}\geq 0$ such that $a_{j}$ is positive for every $j\geq j_{0}$. If $\underline{j}$ is infinite, then the set $\Theta$ is almost surely empty by virtue of Theorem \ref{loidimTheta}. Conversely, let us suppose that $\underline{j}$ is finite and let
\[
d_{*}=\liminf_{j\to\infty}\frac{\log a_{\underline{j}}+\ldots+\log a_{j}}{\log c_{\underline{j}}+\ldots+\log c_{j}}\geq 0.
\]
Theorem \ref{loidimTheta} then ensures that with probability one, the set $\Theta$ either is empty or has Hausdorff measure $d_{*}$. Moreover, $f_{\underline{j}}$ is clearly equal to zero, so that Proposition \ref{altthmprobThetavide} may be applied. This result enables to show that $\Theta$ is empty with probability zero if $b_{j}$ is positive for some $j\geq\underline{j}$ and with probability $\Phi_{\underline{j}}(0)$ otherwise.

\subsubsection{Binary case in dimension one}

Let us consider the particular case in which $d$ is equal to one and the sequence $(c_{j})_{j\geq 0}$ is constant and equal to two. Accordingly, $\etiq_{0}$ is the binary tree $\tree_{2}$ and the compacts $J_{u}$ defined by (\ref{defJuperco}) are the closed dyadic subintervals of $[0,1]$. Here, no specific assumption is made on the form of the transition probability measures $\nu_{t,j}$. Our aim is to illustrate Remark \ref{remarksigma} above, according to which Propositions \ref{condnecsufThetavide1} and \ref{condnecsufThetavide2} may sometimes lead to a criterion to know if the set $\Theta$ is empty with positive probability or with probability one, respectively.

To this end, for any integer $j\geq 0$, let
\[
\left\{\begin{array}{rcl}
\eta_{j} &=& 1-\nu_{0,j}(\{(0,0)\}) \\[2mm]
\gamma_{j} &=& 2\,\nu_{1,j}(\{(1,1)\})+\nu_{1,j}(\{(1,0),(0,1)\}) \\[2mm]
\varsigma_{j} &=& 2\sum\limits_{n=j}^\infty \frac{\nu_{1,n}(\{(1,1)\})}{\gamma_{n}\cdot\prod\limits_{\ell=j}^n\gamma_{\ell}}
\end{array}\right.
\]
and let $\underline{j}$ denote the infimum of all integers $j_{0}\geq 0$ such that $\gamma_{j}>0$ for any $j\geq j_{0}$. Furthermore, let $d_{*}=-\infty$ if $\underline{j}$ is infinite and let
\[
d_{*}=\liminf_{j\to\infty}\frac{\log\gamma_{\underline{j}}+\ldots+\log\gamma_{j}}{j\log 2}
\]
otherwise. In order to illustrate Remark \ref{remarksigma}, let us now establish the following result. Recall that the generating functions $\Phi_{j}$ are defined by (\ref{defPhij}).

\begin{prp}\label{corremarksigma}
If $d_{*}<0$, then $\Theta$ is empty with probability one. If not, then with probability one, $\Theta$ is empty or has Hausdorff dimension $d_{*}$ and, in addition,
\begin{enumerate}
\item\label{corremarksigma1} if there is a $j_{*}\geq 0$ such that $\nu_{1,j}(\{(0,0)\})=0$ for any $j\geq j_{*}$, then
\[
\prob(\Theta=\emptyset)=\Phi_{j_{*}}(0)\cdot\prod_{j=j_{*}}^\infty (1-\eta_{j})^{2^j};
\]
\item\label{corremarksigma2} if $\nu_{1,j}(\{(0,0)\})>0$ for infinitely many integers $j\geq 0$ and if $\sum_{j}2^j\eta_{j}<\infty$, then $\prob(\Theta=\emptyset)$ is positive and it is equal to one if and only if
\[
\varsigma_{\underline{j}}=\infty
\qquad\text{or}\qquad
\liminf_{j\to\infty}\prod_{\ell=\underline{j}}^j \gamma_{\ell}=0
\qquad\text{or}\qquad
\left\{\begin{array}{l}
\Phi_{\underline{j}}(0)=1 \\[2mm]
\forall j\geq\underline{j} \quad \eta_{j}=0;
\end{array}\right.
\]
\item\label{corremarksigma3} if $\nu_{1,j}(\{(0,0)\})>0$ for infinitely many integers $j\geq 0$, if $\sum_{j}2^j\eta_{j}=\infty$ and if $d_{*}>0$, then $\prob(\Theta=\emptyset)$ is less than one and it is equal to zero if and only if $\sum_{j} 2^j\eta_{j}/\varsigma_{j+1}=\infty$.
\end{enumerate}
\end{prp}

\begin{proof}
Theorem \ref{loidimTheta} directly ensures that $\Theta$ is empty with probability one if $d_{*}$ is negative and that, with probability one, $\Theta$ is empty or has Hausdorff dimension $d_{*}$ if $d_{*}$ is nonnegative. Moreover, (\ref{corremarksigma1}) follows at once from Proposition \ref{altthmprobThetavide}. Thus, we may restrict our attention to proving (\ref{corremarksigma2}) and (\ref{corremarksigma3}). To this end, let us assume that $d_{*}\geq 0$ and that $\nu_{1,j}(\{(0,0)\})>0$ for infinitely many integers $j\geq 0$. Then, (\ref{condfjpos}) holds, so that we may apply Propositions \ref{condnecsufThetavide1} and \ref{condnecsufThetavide2} in what follows. In addition, recall that, for any integer $j\geq 0$, the numbers $\underline{\sigma}_{j}$ and $\overline{\sigma}_{j}$ are defined by (\ref{defsigmaj}) and observe that for all $j\geq 0$,
\begin{equation}\label{sigmajencadr}
\left\{\begin{array}{l}
\underline{\sigma}_{j}=\max\Biggl(1,\frac{\varsigma_{j}}{2}+\limsup\limits_{n\to\infty}\frac{1}{\prod\limits_{\ell=j}^n \gamma_{\ell}}\Biggr) \\[2mm]
\overline{\sigma}_{j}=\varsigma_{j}+\liminf\limits_{n\to\infty}\frac{1}{\prod\limits_{\ell=j}^n \gamma_{\ell}}=1+\sum\limits_{n=j}^\infty \frac{2\,\nu_{1,n}(\{(1,1)\})+\gamma_{n}-{\gamma_{n}}^2}{\gamma_{n}\cdot\prod\limits_{\ell=j}^n\gamma_{\ell}}\geq 1.
\end{array}\right.
\end{equation}
In the last sum, $2\,\nu_{1,n}(\{(1,1)\})+\gamma_{n}-{\gamma_{n}}^2$ is nonnegative for any integer $n$, because it is the variance of $Y_{1}+Y_{2}$, when $(Y_{1},Y_{2})$ is distributed according to $\nu_{1,n}$. Note that in particular, if $d_{*}$ is positive, then $\underline{\sigma}_{j}=\max(1,\varsigma_{j}/2)$ and $\overline{\sigma}_{j}=\varsigma_{j}\geq 1$.

Thanks to Proposition \ref{condnecsufThetavide1} and the observation that $\ph_{0,j}(1-z)\geq 1-2\eta_{j}z$ for any real $z\geq 0$ and any integer $j\geq 0$, it is straightforward to check that
\[
\sum_{j=0}^\infty \frac{2^j\eta_{j}}{\frac{\varsigma_{j+1}}{2}+\limsup\limits_{n\to\infty}\frac{1}{\prod\limits_{\ell=j+1}^n \gamma_{\ell}}}<\infty \qquad\Longrightarrow\qquad \prob(\Theta=\emptyset)>0.
\]
Together with (\ref{sigmajencadr}), this implies that $\Theta$ is empty with positive probability in the case where $\sum_{j} 2^j\eta_{j}$ converges and in the case where $\sum_{j} 2^j\eta_{j}/\varsigma_{j+1}$ converges and $d_{*}$ is positive.

Conversely, if $\sum_{j} 2^j\eta_{j}/\varsigma_{j+1}$ diverges while $d_{*}$ remains positive, then $\Theta$ is empty with probability zero, because of (\ref{sigmajencadr}), Proposition \ref{condnecsufThetavide1} and the fact that $\ph_{0,j}(1-z)\leq 1-\eta_{j}z$ for any real $z\in [0,1]$ and any integer $j\geq 0$.

Furthermore, if $\sum_{j} 2^j\eta_{j}=\infty$ and $d_{*}>0$, then $\Theta$ cannot be empty with probability one, in view of (\ref{sigmajencadr}), Proposition \ref{condnecsufThetavide2} and the fact that $\varsigma_{\underline{j}}$ is necessarily finite when $d_{*}$ is positive, because $\nu_{1,n}(\{(1,1)\})\leq\gamma_{n}$ for any integer $n\geq 0$.

Lastly, if $\sum_{j} 2^j\eta_{j}<\infty$, then the criterion to know if $\Theta$ is empty with probability one directly follows from (\ref{sigmajencadr}) and Proposition \ref{condnecsufThetavide2}.
\end{proof}

\begin{rem}
Proposition \ref{corremarksigma} is employed in \cite{Durand:2007kx}, where the set $\Theta$ comes into play in the study of the pointwise regularity of the trajectories of a stochastic process of a certain form. This process is actually a random wavelet series and the correlations between wavelet coefficients are governed by a tree-indexed Markov chain. In that context, the set $\Theta$ is related with the set of points at which the regularity of the process (measured using the notion of H\"older exponent) is the worst possible. Thus, the determination of the Hausdorff dimension of $\Theta$ is crucial for the achievement of the multifractal analysis of the trajectories.
\end{rem}

\subsubsection{Further discussion}

Other generalizations of Mandelbrot's fractal percolation process may be studied with the help of the results of Section \ref{arbrespres}. For instance, one could define a generalization of the random Sierpinski carpet featured in \cite{Dekking:1990lr} and then examine the distribution of its Hausdorff dimension. One could also consider situations where the compacts $J_{u}$ are not simply cubes, but have a more complicated (deterministic or random) geometric structure.

\subsection{A case in which the compacts $J_{u}$ are random}

Until the end of this section, we restrict our attention to an elementary example in which the compacts $J_{u}$ arising in the construction of the set $\Theta$ are random. Many other examples may easily be obtained by generalizing those exposed in \cite{Falconer:1986le,Graf:1988lr,Lyons:2005fc,Mauldin:1986ei}.

Let $p\in (0,1)$ and, for any integer $j\geq 0$, let $\lambda_{j}$ denote a probability measure. We assume that the supports of the measures $\lambda_{j}$ are included in a common proper subinterval of $[0,1]$. The compact subsets $J_{u}$ are now indexed by the binary tree $\tree_{2}$ and are defined as follows. To begin with, the compact $J_{\raci}$ is equal to $[0,1]$. Then, given a random variable $Y_{\raci}$ with distribution $\lambda_{0}$, the compacts $J_{1}$ and $J_{2}$ are equal to $[0,Y_{\raci}]$ and $[Y_{\raci},1]$ respectively. Next, given two random variables $Y_{1}$ and $Y_{2}$ with common distribution $\lambda_{1}$, the compacts $J_{11}$, $J_{12}$, $J_{21}$ and $J_{22}$ are equal to $[0,Y_{\raci}Y_{1}]$, $[Y_{\raci}Y_{1},Y_{\raci}]$, $[Y_{\raci},Y_{\raci}+(1-Y_{\raci})Y_{2}]$ and $[Y_{\raci}+(1-Y_{\raci})Y_{2},1]$ respectively. This procedure is then iterated ad infinitum. It is easy to check that Assertions (\ref{condcompact1})-(\ref{condcompact3}) hold, if for any integer $j\geq 0$, the measure $\mu_{j}$ is the law of $(Y,1-Y)$ when $Y$ is distributed according to $\lambda_{j}$.

Subsequently, some of the compacts $J_{u}$ are either retained or discarded according to the following recursive procedure. To begin with, the compact $J_{\raci}$ is always kept. Next, for any vertex $u\in\tree_{2}$, if the compact $J_{u}$ has been retained, then the compacts $J_{u1}$ and $J_{u2}$ are independently kept, with probability $p$, or discarded, with probability $1-p$, and if the compact $J_{u}$ has been discarded, then the compacts $J_{u1}$ and $J_{u2}$ are thrown away as well. Let us examine the distribution of the Hausdorff dimension of the set of points resulting from the remaining compacts. More formally, this set coincides with the set $\Theta$ defined by (\ref{defTheta}) when the transition probability measures of the Markov chain $(X_{u})_{u\in\tree_{2}}$ are given by
\[
\nu_{0,j}=\delta_{(0,0)} \qquad\text{and}\qquad \nu_{1,j}=(p\,\delta_{1}+(1-p)\delta_{0})^{\otimes 2}
\]
and its initial state is $1$, which means that $X_{\raci}=1$ with probability one. The results of Section \ref{arbrespres} show that if $p\leq 1/2$, then the set $\Theta$ is almost surely empty and that if $p>1/2$, then either the set $\Theta$ is empty or its Hausdorff dimension is the unique solution of
\[
\liminf_{j\to\infty} \frac{1}{j}\sum_{n=0}^{j-1}\log\left(\int_{0}^1\left(y^s+(1-y)^s\right)\lambda_{n}(\dd y)\right)=\log\frac{1}{p}, \qquad s\in\R.
\]
Moreover, in this last case, the set $\Theta$ is empty with probability $(-1+1/p)^2$.

\section{Ancillary results}\label{arbresprel}

In this section, we establish a few lemmas that are called upon at various points of the rest of the paper.

\subsection{Branching processes in varying environment}\label{lienGW}

We first introduce a family of branching processes in varying environment related to the transition probability measures $\nu_{1,j}$ of the Markov chain $(X_{u})_{u\in\etiq_{0}}$ and then establish a relationship between these processes and the set $\Theta$ given by (\ref{defTheta}). Such processes are defined in the same way as the usual Galton-Watson branching processes, except that the offspring distribution of the individuals may depend on their generation, see \cite{Jagers:1974gp}.

\medskip

Recall that the functions $\ph_{1,j}$ are defined by (\ref{defphtj}). Let $(\calX_{j,n,k})_{j,n\geq 0,k\geq 1}$ denote a family of independent random variables such that, for any $j$ and $n$ fixed, $\calX_{j,n,1},\calX_{j,n,2},\ldots$ have generating function $\ph_{1,j+n}$, that is,
\[
\forall k\geq 1 \quad \forall z \qquad \esp[z^{\calX_{j,n,k}}]=\ph_{1,j+n}(z).
\]
Note that with probability one, $\calX_{j,n,k}$ is at most $m_{j+n}$. For any integer $j\geq 0$, the branching process in varying environment $(\calZ_{j,n})_{n\geq 0}$ with offspring distributions having generating functions $\ph_{1,j},\ph_{1,j+1},\ldots$ is then defined by $\calZ_{j,0}=1$ and
\[
\forall n\geq 0 \qquad \calZ_{j,n+1}=\sum_{k=1}^{\calZ_{j,n}} \calX_{j,n,k}.
\]
It is easy to prove by induction on $n$ that the generating function and the expectation of $\calZ_{j,n}$ are respectively given by
\begin{equation}\label{fctgeneespcalZ}
\esp[z^{\calZ_{j,n}}]=\ph_{1,j}\circ\ldots\circ\ph_{1,j+n-1}(z) \qquad\text{and}\qquad \esp[\calZ_{j,n}]=\prod\limits_{\ell=0}^{n-1} \ph_{1,j+\ell}'(1).
\end{equation}

As we shall show in Sections \ref{arbresmajor} and \ref{arbresminor}, the size properties of the set $\Theta$ are closely related to the asymptotic behavior of the processes $(\calZ_{j,n})_{n\geq 0}$. If $j$ is less than the number $\underline{j}$ defined by (\ref{defunderj}), the expectation of $\calZ_{j,n}$ clearly vanishes for all $n$ large enough, so that the process $(\calZ_{j,n})_{n\geq 0}$ becomes extinct (i.e.~$\calZ_{j,n}$ vanishes for all $n$ large enough) with probability one and the study of its asymptotic behavior is elementary. With a view to examining the asymptotic behavior of $(\calZ_{j,n})_{n\geq 0}$ when $j\geq\underline{j}$, let us consider its normed process defined by
\[
\forall n\geq 0 \qquad \calW_{j,n}=\frac{\calZ_{j,n}}{\esp[\calZ_{j,n}]}.
\]
It is straightforward to check that $(\calW_{j,n})_{n\geq 0}$ is a nonnegative martingale. Doob's convergence theorem then ensures that $\calW_{j,n}$ converges almost surely as $n\to\infty$ to a nonnegative random variable denoted by $\calW_{j,\infty}$, see \cite[p. 450]{Doob:1984oi}.

Before going into detail on the asymptotic behavior of $(\calZ_{j,n})_{n\geq 0}$, let us state a consequence of the assumptions made on the random compacts $(J_{u})_{u\in\etiq_{0}}$ coming into play in the construction of the random set $\Theta$.

\begin{lem}\label{lemnjborne}
The sequence $(m_{j})_{j\geq 0}$ is bounded.
\end{lem}

\begin{proof}
Owing to Assertions (\ref{condcompact1}) and (\ref{condcompact2}), there exists a real $C>0$ such that $\sum_{k=1}^{m_{\gene{u}}} \diam{J_{uk}}^d\leq C\diam{J_{u}}^d$ for any vertex $u\in\etiq_{0}$. Assertion (\ref{condcompact3}) then implies that $m_{j}\underline{\beta}^d\leq C$ for any $j\geq 0$. The result follows from the fact that $\underline{\beta}$ is positive.
\end{proof}

Lemma \ref{lemnjborne} is elementary, but crucial, as it ensures that the offspring distributions of the processes $(\calZ_{j,n})_{n\geq 0}$ are uniformly bounded. A result of R. Lyons \cite[Theorem 4.14]{Lyons:1992mz} then immediately implies that $\calW_{j,\infty}$ is positive with probability one given nonextinction. Note that this need not hold for a branching process in varying environment whose offspring distributions are not uniformly bounded, as shown by the various examples given in \cite{DSouza:1994yx,DSouza:1992ms,Lindvall:1974lr,MacPhee:1983fk}.

We use the result of R. Lyons in order to establish the following lemma, which provides all the properties concerning the asymptotic behavior of the processes $(\calZ_{j,n})_{n\geq 0}$ that we shall need in the rest of the paper. Recall that the real $f_{j}$ is defined by (\ref{deffj}) for every integer $j\geq 0$.

\begin{lem}\label{lemlienfj}
For any integer $j\geq 0$,
\[
\left\{\begin{array}{rcl}
j<\underline{j} & \Longrightarrow & \prob(\calZ_{j,n}\to 0 \text{ as } n\to\infty)=1 \\[2mm]
j\geq\underline{j} & \Longrightarrow & \prob(\calZ_{j,n}\to 0 \text{ as } n\to\infty)=\prob(\calW_{j,\infty}=0)=f_{j}.
\end{array}\right.
\]
\end{lem}

\begin{proof}
The expressions of the probability that $\calZ_{j,n}\to 0$ follow from (\ref{fctgeneespcalZ}) and the fact that $\calZ_{j,n}\to 0$ if and only if $\calZ_{j,n}=0$ for all $n$ large enough.

Let us assume that $j\geq\underline{j}$. If $\calZ_{j,n}\to 0$ with probability one, then $\calW_{j,\infty}$ clearly vanishes with probability one. Otherwise, Lemma \ref{lemnjborne} above and Theorem 4.14 in \cite{Lyons:1992mz} imply that $\calW_{j,\infty}$ is positive with probability one given the fact that $\calZ_{j,n}$ does not tend to zero. The result is thus a direct consequence of the observation that $\calW_{j,\infty}$ vanishes if $\calZ_{j,n}$ tends to zero.
\end{proof}

The following result gives some useful properties concerning the reals $f_{j}$. Recall that $\underline{\sigma}_{j}$ and $\overline{\sigma}_{j}$ are defined by (\ref{defsigmaj}).

\begin{lem}\label{lemrecencfj}
For any integer $j\geq 0$,
\[
f_{j}=\ph_{1,j}(f_{j+1}) \qquad\text{and}\qquad 1-\frac{1}{\underline{\sigma}_{j}}\leq f_{j} \leq 1-\frac{1}{\overline{\sigma}_{j}}.
\]
\end{lem}

\begin{proof}
The recurrence relation results from the continuity of $\ph_{1,j}$ and the definition of $f_{j}$. The lower bound on $f_{j}$ is a straightforward consequence of Theorem 2.1 and Proposition 3.1 in \cite{Fujimagari:1980fk}. The upper bound is given by Theorem 1 in \cite{Agresti:1975xc}.
\end{proof}

\medskip

Let us now supply a connection between the set $\Theta$ defined by (\ref{defTheta}) and the branching processes in varying environment $(\calZ_{j,n})_{n\geq 0}$. By definition, $\Theta$ can be written as the union over all vertices $u\in\etiq_{0}$ of the compacts $K_{u}$ given by (\ref{defKu}). As a consequence, the set $\Theta$ is empty if and only if all the compacts $K_{u}$ are empty. Hence, the study of the emptiness probability of $\Theta$ reduces to investigating the probability that all the compacts $K_{u}$ are empty. As shown by the following lemma, this amounts to analyzing the extinction probabilities of the processes $(Z_{u,j})_{j\geq\gene{u}}$ defined by
\begin{equation}\label{defZuj}
\forall j\geq\gene{u} \qquad Z_{u,j}=\#\left\{ v\in\tau_{u} \:\bigl|\: \gene{v}=j \right\}.
\end{equation}

\begin{lem}\label{equiKuvideZujext}
For any vertex $u\in\etiq_{0}$, the compact $K_{u}$ is empty if and only if the process $(Z_{u,j})_{j\geq\gene{u}}$ becomes extinct.
\end{lem}

\begin{proof}
By virtue of (\ref{defKu}), the set $K_{u}$ is empty if and only if the set $\partial\tau_{u}$ is empty. Moreover, the process $(Z_{u,j})_{j\geq\gene{u}}$ becomes extinct if and only if the set $\tau_{u}$ is finite, which is equivalent to the emptiness of $\partial\tau_{u}$. This is due to K\"onig's lemma, according to which a tree in which every vertex has a finite number of sons is finite if and only if its boundary is empty, see \cite{Kleene:2002gf}.
\end{proof}

The next result links the processes $(Z_{u,j})_{j\geq\gene{u}}$ with the branching processes in varying environment $(\calZ_{j,n})_{n\geq 0}$.

\begin{lem}\label{prplienGW}
Let us consider a vertex $u$ in $\etiq_{0}$.
\begin{enumerate}
\item If $\prob(X_{u}=0)>0$, then conditionally on the event $\{X_{u}=0\}$, the process $(Z_{u,j})_{j\geq\gene{u}}$ is almost surely equal to zero.
\item If $\prob(X_{u}=1)>0$, then conditionally on the event $\{X_{u}=1\}$, the process $(Z_{u,j})_{j\geq\gene{u}}$ has the same distribution as the process $(\calZ_{\gene{u},j-\gene{u}})_{j\geq\gene{u}}$.
\end{enumerate}
\end{lem}

\begin{proof}
The first part of the lemma is immediate. In order to establish the second part, observe that $(Z_{u,j})_{j\geq\gene{u}}$ is an inhomogeneous Markov chain with state space $\{0,1,\ldots\}$ such that
\begin{equation}\label{ZujMarkovChain}
\forall j\geq\gene{u} \quad \forall z \qquad \esp[z^{Z_{u,j+1}} \:|\: Z_{u,\gene{u}},\ldots,Z_{u,j}]=\ph_{1,j}(z)^{Z_{u,j}}.
\end{equation}
Indeed, a vertex $v\in\etiq_{0}$ with generation $j+1$ belongs to $\tau_{u}$ if and only if its father $\pi(v)$ belongs to $\tau_{u}$ and $X_{v}=1$. Thus,
\[
Z_{u,j+1}=\sum_{w\in\tau_{u} \atop \gene{w}=j}\sum_{k=1}^{m_{j}} X_{wk}.
\]
The Markov condition (\ref{condMarkov}) ensures that, conditionally on the $\sigma$-field $\trib_{j}$ generated by the variables $X_{w}$ for $w\in\etiq_{0}$ with generation at most $j$, the vectors $(X_{wk})_{k\in\intn{1}{m_{j}}}$ for $w\in\tau_{u}$ with generation $j$ are independent and distributed according to $\nu_{1,j}$. Therefore, 
\begin{equation}\label{Zuj1sum}
\esp[z^{Z_{u,j+1}} \:|\: \trib_{j}]=\ph_{1,j}(z)^{Z_{u,j}}.
\end{equation}
The tower property of conditional expectation and the fact that the variables $Z_{u,\gene{u}},\ldots,Z_{u,j}$ are $\trib_{j}$-measurable then lead to (\ref{ZujMarkovChain}). The second part of the lemma follows from the observation that (\ref{ZujMarkovChain}) also holds for the process $(\calZ_{\gene{u},j-\gene{u}})_{j\geq\gene{u}}$, that $\calZ_{\gene{u},0}=1$ and that $Z_{u,\gene{u}}=1$ if and only if $X_{u}=1$.
\end{proof}

\subsection{Flows in random networks}\label{lienecoul}

K. Falconer \cite{Falconer:1986le} observed that the problem of determining the $s$-dimensional Hausdorff measures of the sets obtained through certain random recursive constructions can be reduced to that of examining flows in random networks. His approach can actually be adapted for studying the Hausdorff dimension of the set $\Theta$ defined by (\ref{defTheta}). Note that we do not go into detail about network theory here. We refer to \cite{Bollobas:1998lr,Ford:1962wt,Lyons:2005fc} for full expositions of this topic.

As $\Theta$ is the union over all vertices $u\in\etiq_{0}$ of the sets $K_{u}$ given by (\ref{defKu}), we shall begin by studying the Hausdorff dimension of $K_{u}$ for any $u\in\etiq_{0}$. Recall that the set $K_{u}$ is based on the boundary $\partial\tau_{u}$ of the set $\tau_{u}$ defined by (\ref{deftauu}). Let $\coup(\tau_{u})$ denote the collection of all finite sets $\chi\subseteq\tau_{u}$ enjoying
\[
\left\{ \begin{array}{lll}
\forall\zeta\in\partial\tau_{u} & \exists v\in\chi & \quad \zeta_{1}\ldots\zeta_{\gene{v}}=v \\
\forall v\in\chi & \exists \zeta\in\partial\tau_{u} & \quad \zeta_{1}\ldots\zeta_{\gene{v}}=v \\
\forall v\in\chi & \forall \ell\in\intn{\gene{u}}{\gene{v}-1} & \quad v_{1}\ldots v_{\ell}\notin\chi.
\end{array}\right.
\]
According to the terminology of network theory, any element of $\coup(\tau_{u})$ is called a {\em cut} of $\tau_{u}$. One easily checks that $\coup(\tau_{u})=\{ \emptyset \}$ if $\partial\tau_{u}$ is empty and that the singleton $\{ u \}$ necessarily belongs to $\coup(\tau_{u})$ otherwise. In addition, for any integer $j\geq\gene{u}$, let $\coup_{j}(\tau_{u})$ denote the collection of all cuts $\chi\in\coup(\tau_{u})$ formed by vertices with generation at least $j$ only.

For every positive real $s$, let us consider
\begin{equation}\label{defEsu}\begin{split}
E_{s,u} &=\inf_{\chi\in\coup(\tau_{u})} \sum_{v\in\chi} \left( \prod_{\ell=\gene{u}+1}^{\gene{v}} L_{v_{1}\ldots v_{\ell}} \right)^s \\[2mm]
\text{and} \quad \tilde E_{s,u} &=\lim_{j\uparrow\infty}\uparrow \inf_{\chi\in\coup_{j}(\tau_{u})} \sum_{v\in\chi} \left( \prod_{\ell=\gene{u}+1}^{\gene{v}} L_{v_{1}\ldots v_{\ell}} \right)^s.
\end{split}\end{equation}
Then, $E_{s,u}$ can actually be seen as a maximal flow through a certain random network associated with the set $\tau_{u}$. The following result shows that, together with $\tilde E_{s,u}$, it is linked with the $s$-dimensional Hausdorff measure of the set $K_{u}$.

\begin{lem}\label{prplienEusHsKu}
There exists a real $C>0$ such that with probability one, for every vertex $u\in\etiq_{0}$ and every real $s>0$,
\[
C\diam{J_{u}}^s E_{s,u}\leq C\diam{J_{u}}^s\tilde E_{s,u}\leq\hau^s(K_{u})\leq \diam{J_{u}}^s\tilde E_{s,u}.
\]
\end{lem}

\begin{proof}
Assertion (\ref{condcompact3}) ensures that with probability one,
\begin{equation}\label{rapdiambeta}
\forall v\in\etiq_{0}^* \qquad \underline{\beta}\leq L_{v}\leq\overline{\beta}.
\end{equation}
Throughout the proof, we assume that the event on which (\ref{rapdiambeta}) holds occurs. Let $\eps>0$ and $j\geq 0$ with $\diam{J_{\raci}}\overline{\beta}^j<\eps$ and let $\chi\in\coup_{j}(\tau_{u})$. Thanks to (\ref{defxzeta}), (\ref{defKu}) and (\ref{rapdiambeta}), it is straightforward to check that the compact sets $J_{v}$, for $v\in\chi$, cover $K_{u}$ and have diameter less than $\eps$. Hence,
\[
\hau^s_{\eps}(K_{u})\leq\sum_{v\in\chi} \diam{J_{v}}^s = \diam{J_{u}}^s \sum_{v\in\chi} \left( \prod_{\ell=\gene{u}+1}^{\gene{v}} L_{v_{1}\ldots v_{\ell}} \right)^s.
\]
We get $\hau^s(K_{u})\leq\diam{J_{u}}^s\tilde E_{s,u}$ by taking the infimum over all $\chi\in\coup_{j}(\tau_{u})$ and letting $j\to\infty$ in the right-hand side and then by letting $\eps\to 0$ in the left-hand side.

Conversely, let $j\geq\gene{u}$ and $\eps\in (0,\diam{J_{\raci}}\underline{\beta}^j]$ and let $(U_{p})_{p\geq 1}$ denote a cover of $K_{u}$ by sets of diameter less than $\eps$. For any $p$, if $\diam{U_{p}}>0$, then let $\tilde U_{p}=U_{p}$, else let $\tilde U_{p}$ denote a set with diameter $\eps/2^p$ that contains $U_{p}$. Thus, $(\tilde U_{p})_{p\geq 1}$ is a cover of $K_{u}$ enjoying $0<\diam{\tilde U_{p}}<\eps$ for all $p$. Moreover, for any $p\geq 1$, let
\[
V_{p}=\bigl\{ v\in\tau_{u}\setminus\{\raci\} \:\bigl|\: \tilde U_{p}\cap J_{v}\neq\emptyset\text{ and }\diam{J_{v}}\leq\diam{\tilde U_{p}}<\diam{J_{\pi(v)}}\bigr\}.
\]
Let us show that $\# V_{p}$ is bounded. Let $x_{0}\in\tilde U_{p}$ and let $\kappa'$ denote a positive real such that $\|x\|_{\infty}\leq\kappa'\|x\|$ for all $x\in\R^d$, where $\|\cdot\|_{\infty}$ and $\|\cdot\|$ are respectively the supremum norm and the norm $\R^d$ is endowed with. The sets $J_{v}$, for $v\in V_{p}$, have disjoint interiors and are included in the closed ball $B$ with center $x_{0}$ and radius $2\kappa'\diam{\tilde U_{p}}$ in the sense of the supremum norm, so that
\[
(4\kappa')^d\diam{\tilde U_{p}}^d=\leb^d(B)\geq\sum_{v\in V_{p}} \leb^d(\interieur{J_{v}})\geq\kappa\sum_{v\in V_{p}} \diam{J_{v}}^d,
\]
where $\interieur{}$ denotes interior and $\kappa$ is given by Assertion (\ref{condcompact2}). As $\diam{J_{v}}\geq\underline{\beta}\diam{J_{\pi(v)}}>\underline{\beta}\diam{\tilde U_{p}}$ for any $v\in V_{p}$, it follows that
\begin{equation}\label{lemgeomeq2}
\# V_{p}\leq \frac{1}{C} \qquad\text{with}\qquad C=\frac{\kappa\underline{\beta}^d}{(4\kappa')^d}.
\end{equation}

Let $\chi$ denote the set obtained by removing from $\chi'=\bigcup_{p} V_{p}$ the vertices $v$ such that $v_{1}\ldots v_{\ell}\in\chi'$ for some $\ell<\gene{v}$ or such that $\zeta_{1}\ldots\zeta_{\gene{v}}\neq v$ for all $\zeta\in\partial\tau_{u}$. One can straightforwardly check that $\chi$ is a cut of $\tau_{u}$. Moreover, $\diam{J_{v}}<\eps\leq\diam{J_{\raci}}\underline{\beta}^j$ for any $v\in\chi$, so that $\underline{\beta}^{\gene{v}}\leq\underline{\beta}^j$ by (\ref{rapdiambeta}). As a result, $\chi$ actually belongs to $\coup_{j}(\tau_{u})$. Furthermore, thanks to (\ref{lemgeomeq2}),
\begin{equation*}\begin{split}
\diam{J_{u}}^s \sum_{v\in\chi} \left( \prod_{\ell=\gene{u}+1}^{\gene{v}} L_{v_{1}\ldots v_{\ell}} \right)^s &\leq \sum_{p=1}^\infty\sum_{v\in V_{p}} \diam{J_{v}}^s \\
&\leq \frac{1}{C}\sum_{p=1}^\infty \diam{\tilde U_{p}}^s\leq \frac{1}{C}\left( \frac{\eps^s}{2^s-1} + \sum_{p=1}^\infty \diam{U_{p}}^s \right).
\end{split}\end{equation*}
We finally obtain $C\diam{J_{u}}^s\tilde E_{s,u}\leq\hau^s(K_{u})$ by taking the infimum over $\chi\in\coup_{j}(\tau_{u})$ in the left-hand side and the infimum over $(U_{p})_{p\geq 1}$ in the right-hand side and by letting $\eps\to 0$ and $j\to\infty$.

To end the proof of the lemma, it suffices to observe that $\coup_{j}(\tau_{u})\subseteq\coup(\tau_{u})$ for any $j\geq\gene{u}$, so that $E_{s,u}\leq\tilde E_{s,u}$.
\end{proof}

\begin{rem}\label{remarkmeasurability}
Note that both $E_{s,u}$ and $\tilde E_{s,u}$ are measurable with respect to the $\sigma$-field generated by the ratios $(L_{v})_{v\in\etiq_{0}^*}$ and the Markov chain $(X_{v})_{v\in\etiq_{0}}$. Lemma \ref{prplienEusHsKu} then implies that the Hausdorff dimension of each $K_{u}$ is measurable with respect to this $\sigma$-field. As $\Theta$ is the union over all vertices $u\in\etiq_{0}$ of the sets $K_{u}$, its Hausdorff dimension is also measurable with respect to the same $\sigma$-field. As observed by K. Falconer \cite{Falconer:1986le}, as far as dimension calculations are concerned, it is not necessary to impose measurability conditions on the compacts $J_{v}$ themselves.
\end{rem}

Lemma \ref{prplienEusHsKu} reduces the problem of computing the Hausdorff dimension of $K_{u}$ to that of determining for which values of $s$ the random variables $E_{s,u}$ and $\tilde E_{s,u}$ are positive. In particular, if the flow $E_{s,u}$ is positive for some $s>0$, then the $s$-dimensional Hausdorff measure of $K_{u}$ is positive as well, so that the Hausdorff dimension of this set is at least $s$. Thus, with a view to later deriving a lower bound on $\dim K_{u}$, we now study the probability that $E_{s,u}$ vanishes.

To this end, for any integer $j\geq 0$, we need to introduce two independent families $(L^j_{u})_{u\in\etiq^*_{j}}$ and $(X^j_{u})_{u\in\etiq_{j}}$ of random variables which satisfy the following conditions, which are similar to those given by Assertions (\ref{condcompact3}) and (\ref{condMarkov}):
\renewcommand{\theenumi}{\Alph{enumi}}
\begin{enumerate}\setcounter{enumi}{4}
\item\label{condcompact4bis} The vectors $(L^j_{uk})_{k\in\intn{1}{m_{j+\gene{u}}}}$, for $u\in\etiq_{j}$, are independent and distributed according to $\mu_{j+\gene{u}}$.
\item\label{condMarkovbis} With probability one, $X^j_{\raci}=1$ and, for any vertex $u\in\etiq_{j}$, the conditional distribution of the vector $(X^j_{uk})_{k\in\intn{1}{m_{j+\gene{u}}}}$, conditionally on $X_{v}^j$ for $v\in\etiq_{j}\setminus(u\etiq^*_{j+\gene{u}})$, is $\nu_{X^j_{u},j+\gene{u}}$.
\end{enumerate}
\renewcommand{\theenumi}{\alph{enumi}}
The last condition means that $(X^j_{u})_{u\in\etiq_{j}}$ is a Markov chain with initial state $1$ and transition probability measures $\nu_{t,j+n}$ for $n\geq 0$ and $t\in\{0,1\}$. For any $u\in\etiq_{j}$, we also need to consider the set
\begin{equation}\label{deftauj0u}
\tau^j_{u}=\bigl\{ v\in u\etiq_{j+\gene{u}} \:|\: \forall n\in\intn{\gene{u}}{\gene{v}} \quad X^j_{v_{1}\ldots v_{n}}=1 \bigr\},
\end{equation}
which is defined as in (\ref{deftauu}) and the collection $\coup(\tau^j_{u})$ of all its cuts. Then, Assertions (\ref{condcompact3}) and (\ref{condMarkov}) imply that for any $u\in\etiq_{0}$ with $\prob(X_{u}=1)>0$ and any $s>0$, conditionally on the event $\{X_{u}=1\}$, the flow $E_{s,u}$ has the same distribution as
\begin{equation}\label{Esuegalenloi}
\inf_{\chi\in\coup(\tau^{\gene{u}}_{\raci})}\sum_{v\in\chi}\left( \prod_{\ell=1}^{\gene{v}} L^{\gene{u}}_{v_{1}\ldots v_{\ell}} \right)^s.
\end{equation}
Moreover, if $X_{u}=0$, then $E_{s,u}$ obviously vanishes. Thus, the problem is reduced to the study of the probabilities
\begin{equation}\label{defesj}
e_{s,j}=\prob\left( \inf_{\chi\in\coup(\tau^j_{\raci})}\sum_{v\in\chi}\left( \prod_{\ell=1}^{\gene{v}} L^j_{v_{1}\ldots v_{\ell}} \right)^s=0 \right),
\end{equation}
for $s>0$ and $j\geq 0$. The following result shows that, for any fixed $s>0$, the reals $e_{s,j}$ satisfy the same recurrence relation as that enjoyed by the reals $f_{j}$, see Lemma~\ref{lemrecencfj}.

\begin{lem}\label{lemrecesj}
For any real $s>0$ and any integer $j\geq 0$,
\[
e_{s,j}=\ph_{1,j}(e_{s,j+1}).
\]
\end{lem}

\begin{proof}
Let $S^j_{1}$ denote the set of vertices $u\in\etiq_{j}$ with $\gene{u}=1$ and $X^j_{u}=1$. If the set $\tau^j_{\raci}$ is finite, then $\coup(\tau^j_{\raci})$ is reduced to the singleton $\{ \emptyset \}$, as well as the sets $\coup(\tau^j_{u})$, for $u\in S^j_{1}$. Conversely, if $\tau^j_{\raci}$ is infinite, then $\coup(\tau^j_{\raci})$ consists of the singleton $\{ \raci \}$, together with all the possible unions of elements of $\coup(\tau^j_{u})$, for $u\in S^j_{1}$. In both cases, it follows that
\[
\inf_{\chi\in\coup(\tau^j_{\raci})}\sum_{v\in\chi}\left( \prod_{\ell=1}^{\gene{v}} L^j_{v_{1}\ldots v_{\ell}} \right)^s=\min\left( 1 , \sum_{u\in S^j_{1}} (L^j_{u})^s \inf_{\chi_{u}\in\coup(\tau^j_{u})} \sum_{v\in\chi_{u}} \left( \prod_{\ell=2}^{\gene{v}} L^j_{v_{1}\ldots v_{\ell}} \right)^s \right).
\]
In particular, the infimum in the left-hand side vanishes if and only if all the infimums in the right-hand side do. Meanwhile, conditionally on the variables $X^j_{u}$ for $\gene{u}\leq 1$, each of these infimums vanishes independently of the others with probability $e_{s,j+1}$. Hence,
\[
\prob\left( \inf_{\chi\in\coup(\tau^j_{\raci})}\sum_{v\in\chi}\left( \prod_{\ell=1}^{\gene{v}} L^j_{v_{1}\ldots v_{\ell}} \right)^s=0 \:\Biggl|\: X^j_{u},\,\gene{u}\leq 1 \right)={e_{s,j+1}}^{\# S^j_{1}}.
\]
In order to conclude, it suffices to observe that $\# S^j_{1}=X^j_{1}+\ldots+X^j_{m_{j}}$ and to take expectations.
\end{proof}

Let us now assume that the number $d_{*}$ defined by (\ref{defdetoile}) is positive. We end this subsection by giving an upper bound on $e_{s,j}$, when $s$ is less than $d_{*}$. For any real $s\in (0,d_{*})$ and any integer $j\geq 0$, let us consider the function $\phi_{s,j}$ defined by
\[
\phi_{s,j}:z\mapsto\int_{\{0,1\}^{m_{j}}}\int_{[\underline{\beta},\overline{\beta}]^{m_{j}}} \prod_{k=1}^{m_{j}} \left( 1-{\ell_{k}}^s(1-z^{x_{k}}) \right) \mu_{j}(\dd\ell) \nu_{1,j}(\dd x).
\]

\begin{lem}\label{lemmajesj}
If $d_{*}>0$, then for any real $s\in (0,d_{*})$ and any integer $j\geq 0$,
\[
e_{s,j}\leq 1-\frac{1}{\varsigma_{s,j}} \qquad\text{where}\qquad \varsigma_{s,j}=\sum_{n=0}^\infty \frac{ \phi_{s,j+n}''(1) }{ \phi_{s,j+n}'(1) \prod\limits_{\ell=0}^n \phi_{s,j+\ell}'(1) }.
\]
\end{lem}

We refer to the next subsection for a proof of this result. This mainly consists in adapting to our setting some techniques of percolation theory exposed by R. Lyons and Y. Peres in \cite[Chapter 4]{Lyons:2005fc}.

\subsection{Percolation on trees} With a view to proving Lemma \ref{lemmajesj}, we begin by considering percolation on the tree $\tau^j_{\raci}$ defined by (\ref{deftauj0u}). Let $\xi$ denote a mapping from $\etiq_{j}^*$ to $\{0,1\}$. This mapping is intended to indicate which vertices of $\tau^j_{\raci}$ remain during the percolation process. Actually, the remaining vertices are those of the set
\[
\xi\cdot\tau^j_{\raci}=\bigl\{ u\in\tau^j_{\raci} \:\bigl|\: \forall j\in\intn{1}{\gene{u}} \quad \xi_{u_{1}\ldots u_{j}}=1 \bigr\}.
\]
This set is the largest subtree of $\tau^j_{\raci}$ formed by the root $\raci$ and the vertices $u\in\tau^j_{\raci}$ for which $\xi_{u}=1$. For any integer $n\geq 0$, let
\[
\xi\cdot Z^j_{\raci,n}=\#\bigl\{ u\in\xi\cdot\tau^j_{\raci} \:\bigl|\: \gene{u}=n \bigr\}.
\]
If the mapping $\xi$ is chosen according to the random product measure
\[
\Mu^j_{s}=\bigotimes_{u\in\etiq_{j}^*} \left( (L^j_{u})^s \delta_{1}+\left( 1-(L^j_{u})^s \right)\delta_{0} \right),
\]
then it is possible to express the generating functions of $\xi\cdot Z^j_{\raci,0},\xi\cdot Z^j_{\raci,1},\ldots$ in terms of the functions $\phi_{s,j},\phi_{s,j+1},\ldots$ This is the purpose of the following result.

\begin{lem}\label{lemlienGWprocperco}
If $d_{*}>0$, then for any real $s\in (0,d_{*})$ and any integer $j\geq 0$,
\[
\forall n\geq 0 \quad \forall z \qquad \esp\left[ \int_{\{0,1\}^{\etiq_{j}^*}} z^{\xi\cdot Z^j_{\raci,n}} \Mu^j_{s}(\dd\xi)\right]=\phi_{s,j}\circ\ldots\circ\phi_{s,j+n-1}(z).
\]
\end{lem}

\begin{proof}
We prove the result by induction on $n\geq 0$. First, the equality is obviously verified for every $z$ when $n=0$. Then, let $n$ denote an integer for which the equality holds for every $z$. For the sake of clarity, we need to introduce some further notations. Let $\Xi_{n}$ be the set of all mappings that are valued in $\{0,1\}$ and defined on the set of vertices $u\in\etiq_{j}^*$ with generation at most $n$ and let us consider the random measure
\[
\Mu^j_{s,n} =\bigotimes_{u\in\etiq_{j}^* \atop \gene{u}\leq n} \left( (L^j_{u})^s \delta_{1}+\left( 1-(L^j_{u})^s \right)\delta_{0} \right).
\]
In addition, let $\tilde\Xi_{n}=\{0,1\}^{m_{j+n}}$ and, for any vertex $w\in\etiq_{j}$ with generation $n$, let us consider the random measure
\[
\tilde \Mu^j_{s,w,n} =\bigotimes_{k=1}^{m_{j+n}} \left( (L^j_{wk})^s \delta_{1}+\left( 1-(L^j_{wk})^s \right)\delta_{0} \right).
\]
For any mapping $\xi:\etiq_{j}^*\to\{0,1\}$, a vertex $v\in\etiq_{j}$ with generation $n+1$ belongs to the tree $\xi\cdot\tau^j_{\raci}$ if and only if its father $\pi(v)$ also belongs to it and if $\xi_{v}=X^j_{v}=1$. As a result,
\[
\xi\cdot Z^j_{\raci,n+1}=\sum_{w\in\xi\cdot\tau^j_{\raci} \atop \gene{w}=n}\sum_{k=1}^{m_{j+n}} \xi_{wk}X^j_{wk}.
\]
Thus, for any complex number $z$, the integral of $z^{\xi\cdot Z^j_{\raci,n+1}}$ with respect to the random product measure $\Mu^j_{s}(\dd\xi)$ is equal to
\begin{equation*}\begin{split}
& \int_{\Xi_{n}} \prod_{w\in\xi\cdot\tau^j_{\raci} \atop \gene{w}=n} \left(\int_{\tilde\Xi_{n}} \prod_{k=1}^{m_{j+n}} \left(z^{X^j_{wk}}\right)^{\xi^w_{k}} \tilde \Mu^j_{s,w,n}(\dd\xi^w) \right) \Mu^j_{s,n}(\dd\xi)\\
=&\int_{\Xi_{n}} \prod_{w\in\xi\cdot\tau^j_{\raci} \atop \gene{w}=n} \prod_{k=1}^{m_{j+n}} \left( (L^j_{wk})^s z^{X^j_{wk}}+1-(L^j_{wk})^s \right) \Mu^j_{s,n}(\dd\xi).
\end{split}\end{equation*}
Therefore, owing to Assertions (\ref{condcompact4bis}) and (\ref{condMarkovbis}), the conditional expectation of the right-hand side of the previous equality, conditionally on the variables $X^j_{u}$ and $L^j_{u}$ for $\gene{u}\leq n$, is equal to
\[
\int_{\Xi_{n}} \phi_{s,j+n}(z)^{\xi\cdot Z^j_{\raci,n}} \Mu^j_{s,n}(\dd\xi)
\]
As a consequence,
\[
\esp\left[ \int_{\{0,1\}^{\etiq_{j}^*}} z^{\xi\cdot Z^j_{\raci,n+1}} \Mu^j_{s}(\dd\xi) \:\Biggl| \begin{array}{l} X^j_{u},\,L^j_{u},\\ \gene{u}\leq n\end{array} \right]=\int_{\{0,1\}^{\etiq_{j}^*}} \phi_{s,j+n}(z)^{\xi\cdot Z^j_{\raci,n}} \Mu^j_{s}(\dd\xi).
\]
It finally suffices to take expectations in order to conclude.
\end{proof}

We are now able to prove Lemma \ref{lemmajesj}. Let us assume that $d_{*}$ is positive, consider a real $s\in (0,d_{*})$ and establish that the probability $e_{s,j}$ is at most $1-1/\varsigma_{s,j}$. We can clearly suppose that $j\geq\underline{j}$, since $\varsigma_{s,j}$ is infinite otherwise.

To begin with, observe that the mean number of vertices remaining in any cut $\chi\in\coup(\tau^j_{\raci})$ after the percolation process has occurred is
\[
\int_{\{0,1\}^{\etiq_{j}^*}} \#(\xi\cdot\tau^j_{\raci}\cap\chi) \Mu^j_{s}(\dd\xi) = \sum_{v\in\chi} \int_{\{0,1\}^{\etiq_{j}^*}} \ind_{\{v\in\xi\cdot\tau^j_{\raci}\}} \Mu^j_{s}(\dd\xi).
\]
Furthermore, any vertex $v\in\chi\subseteq\tau^j_{\raci}$ also belongs to $\xi\cdot\tau^j_{\raci}$ if and only if $\xi_{v_{1}\ldots v_{\ell}}=1$, for every integer $\ell\in\intn{1}{\gene{v}}$. Hence,
\[
\int_{\{0,1\}^{\etiq_{j}^*}} \ind_{\{v\in\xi\cdot\tau^j_{\raci}\}} \Mu^j_{s}(\dd\xi)=\left(\prod_{\ell=1}^{\gene{v}} L^j_{v_{1}\ldots v_{\ell}}\right)^s.
\]
The mean number of vertices remaining in the cut $\chi$ is then obtained by summing the right-hand side over all vertices $v\in\chi$. Meanwhile, this mean number is at least
\[
\int_{\{0,1\}^{\etiq_{j}^*}} \ind_{\{ \#(\xi\cdot\tau^j_{\raci} \cap \chi) \geq 1 \}} \Mu^j_{s}(\dd\xi) \geq \int_{\{0,1\}^{\etiq_{j}^*}} \ind_{\{ \xi\cdot Z^j_{\raci,n} \nrightarrow 0\text{ as }n\to\infty \}} \Mu^j_{s}(\dd\xi).
\]
Indeed, the mean number of vertices remaining in $\chi$ is greater than or equal to the probability that there remains at least one vertex in $\chi$. In addition, if $\xi\cdot Z^j_{\raci,n}$ does not tend to zero as $n\to\infty$, the boundary of the tree $\xi\cdot\tau^j_{\raci}$, which is defined as in (\ref{deffrtauu}), contains at least a sequence $\zeta=(\zeta_{j})_{j\geq 1}$, by virtue of K\"onig's lemma. This sequence also belongs to the boundary of the tree $\tau^j_{\raci}$, so that the cut $\chi$ contains a vertex $v$ enjoying $\zeta_{1}\ldots\zeta_{\gene{v}}=v$. The vertex $v$ thus simultaneously belongs to $\chi$ and $\xi\cdot\tau^j_{\raci}$. Therefore, at least a vertex remains in $\chi$.

Taking the infimum over all cuts $\chi$ in $\coup(\tau^j_{\raci})$, we deduce that
\[
\inf_{\chi\in\coup(\tau^j_{\raci})}\sum_{v\in\chi} \left(\prod_{\ell=1}^{\gene{v}} L^j_{v_{1}\ldots v_{\ell}}\right)^s \geq \int_{\{0,1\}^{\etiq_{j}^*}} \ind_{\{ \xi\cdot Z^j_{\raci,n} \nrightarrow 0\text{ as }n\to\infty \}} \Mu^j_{s}(\dd\xi).
\]
In particular, if this infimum vanishes, then the preceding integral vanishes as well. Owing to the definition (\ref{defesj}) of the probability $e_{s,j}$, this observation implies that
\begin{equation*}\begin{split}
e_{s,j} &\leq \prob\left( \int_{\{0,1\}^{\etiq_{j}^*}} \ind_{\{ \xi\cdot Z^j_{\raci,n} \nrightarrow 0\text{ as }n\to\infty \}} \Mu^j_{s}(\dd\xi)=0 \right) \\
& \leq \esp\left[ \int_{\{0,1\}^{\etiq_{j}^*}} \ind_{\{ \xi\cdot Z^j_{\raci,n} \to 0\text{ as }n\to\infty \}} \Mu^j_{s}(\dd\xi) \right].
\end{split}\end{equation*}
Observe that $\xi\cdot Z^j_{\raci,n}$ tends to zero as $n\to\infty$ if and only if $\xi\cdot Z^j_{\raci,n}=0$ for all $n$ large enough. Owing to Lemma \ref{lemlienGWprocperco}, the last expectation is thus the limit of
\[
\esp\left[ \int_{\{0,1\}^{\etiq_{j}^*}} \ind_{\{ \xi\cdot Z^j_{\raci,j_{1}}=0 \}} \Mu^j_{s}(\dd\xi) \right]=\phi_{s,j}\circ\ldots\circ\phi_{s,j+j_{1}-1}(0)
\]
as $j_{1}\to\infty$. Furthermore, for any $j_{1}\geq 0$, Theorem 1 in \cite{Agresti:1975xc} ensures that the right-hand side is at most
\[
1-\left(\frac{ 1 }{ \prod\limits_{\ell=0}^{j_{1}-1} \phi_{s,j+\ell}'(1) }+\sum_{n=0}^{j_{1}-1} \frac{ \phi_{s,j+n}''(1) }{ \phi_{s,j+n}'(1) \prod\limits_{\ell=0}^n \phi_{s,j+\ell}'(1) }\right)^{-1}.
\]
To end the proof of Lemma \ref{lemmajesj}, it remains to observe that this expression tends to $1-1/\varsigma_{s,j}$ as $j_{1}\to\infty$. This is due to the fact that $s$ is less than $d_{*}$, together with the observation that $\phi_{s,j+\ell}'(1)$ is equal to the number $\alpha_{s,j+\ell}$ defined by (\ref{defalphasj}), for any nonnegative integer $\ell$.

\section{Upper bound on the dimension}\label{arbresmajor}

Recall that, by virtue of its definition (\ref{defTheta}), the set $\Theta$ is the union over all vertices $u\in\etiq_{0}$ of the compacts $K_{u}$ given by (\ref{defKu}). Hence, with a view to proving Theorem \ref{loidimTheta}, we establish in this section that the Hausdorff dimension of the compacts $K_{u}$ is at most the number $d_{*}$ defined by (\ref{defdetoile}). We first discuss the elementary case in which the generation of the vertex $u$ is less than the number $\underline{j}$ defined by (\ref{defunderj}).

\begin{prp}\label{majdimKueasy}
For any vertex $u\in\etiq_{0}$ with $\gene{u}<\underline{j}$,
\[
\as \qquad \dim K_{u}=-\infty.
\]
\end{prp}

\begin{proof}
Lemmas \ref{lemlienfj}, \ref{equiKuvideZujext} and \ref{prplienGW} ensure that for any vertex $u\in\etiq_{0}$ with generation less than $\underline{j}$, the compact $K_{u}$ is almost surely empty. Thus, its Hausdorff dimension is $-\infty$ with probability one.
\end{proof}

Let us now consider the case in which the generation of $u$ is at least $\underline{j}$. The number $\underline{j}$ is thus necessarily finite. Lemma \ref{prplienEusHsKu} ensures that, in order to derive an upper bound on the Hausdorff dimension of $K_{u}$, it suffices to identify values of the positive real $s$ for which the random variable $\tilde E_{s,u}$ vanishes. The next lemma suggests that this may be done by examining the asymptotic behavior of the processes $(Z_{s,u,j})_{j\geq\gene{u}}$ defined by
\begin{equation}\label{defZsuj}
\forall s>0 \qquad \forall j\geq\gene{u} \qquad Z_{s,u,j}=\sum_{v\in\tau_{u} \atop \gene{v}=j} \left( \prod_{\ell=\gene{u}+1}^j L_{v_{1}\ldots v_{\ell}} \right)^s,
\end{equation}
where the set $\tau_{u}$ is defined by (\ref{deftauu}) and the ratios $L_{v_{1}\ldots v_{\ell}}$ are given by (\ref{defLu}).

\begin{lem}\label{tildeEsuleqliminfZsuj}
For any vertex $u\in\etiq_{0}$ with $\gene{u}\geq\underline{j}$ and any real $s>0$,
\[
\tilde E_{s,u}\leq\liminf_{j\to\infty} Z_{s,u,j}.
\]
\end{lem}

\begin{proof}
For each integer $j\geq\gene{u}$, the set of vertices $v\in\tau_{u}$ with generation $j$ for which $Z_{v,n}$ does not tend to zero as $n\to\infty$ belongs to $\coup_{j}(\tau_{u})$, the collection of all cuts of $\tau_{u}$ formed by vertices with generation at least $j$ only. Consequently,
\[
\inf_{\chi\in\coup_{j}(\tau_{u})}\sum_{v\in\chi} \left( \prod_{\ell=\gene{u}+1}^{\gene{v}} L_{v_{1}\ldots v_{\ell}} \right)^s \leq Z_{s,u,j}.
\]
The result is then a straightforward consequence of (\ref{defEsu}).
\end{proof}

In what follows, we also consider the process $(Z_{0,u,j})_{j\geq\gene{u}}$ obtained by letting $s=0$ in (\ref{defZsuj}). It is clearly equal to the process $(Z_{u,j})_{j\geq\gene{u}}$ given by (\ref{defZuj}). In addition, we make use of the normed processes $(W_{s,u,j})_{j\geq\gene{u}}$ defined by
\begin{equation}\label{defWsuj}
\forall s\geq 0 \quad \forall j\geq\gene{u} \qquad W_{s,u,j}=\frac{Z_{s,u,j}}{\prod\limits_{\ell=\gene{u}}^{j-1} \alpha_{s,\ell}},
\end{equation}
where the numbers $\alpha_{s,\ell}$ are given by (\ref{defalphasj}). These numbers are positive if $\ell\geq\underline{j}$, so that the normed processes are correctly defined.

\begin{lem}\label{lemWsujmart}
For any vertex $u\in\etiq_{0}$ with $\gene{u}\geq\underline{j}$ and any real $s\geq 0$, the process $(W_{s,u,j})_{j\geq\gene{u}}$ is a nonnegative martingale.
\end{lem}

\begin{proof}
For any integer $j\geq\gene{u}$, a vertex $w$ in $\etiq_{0}$ with generation $j+1$ belongs to $\tau_{u}$ if and only if its father $\pi(w)$ also belongs to $\tau_{u}$ and if $X_{w}=1$. Hence,
\[
Z_{s,u,j+1}=\sum_{v\in\tau_{u} \atop \gene{v}=j} \left( \prod_{\ell=\gene{u}+1}^j L_{v_{1}\ldots v_{\ell}} \right)^s\,\sum_{k=1}^{m_{j}} {L_{vk}}^s X_{vk}.
\]
Assertions (\ref{condcompact3}) and (\ref{condMarkov}) then imply that the conditional expectation of $Z_{s,u,j+1}$, conditionally on the $\sigma$-field generated by the variables $X_{v}$ and $\diam{J_{v}}$ for $v\in\etiq_{0}$ such that $\gene{v}\leq j$, is equal to $\alpha_{s,j} Z_{s,u,j}$. The result follows from the fact that the variables $W_{s,u,\gene{u}},\ldots,W_{s,u,j}$ are measurable with respect to this $\sigma$-field.
\end{proof}

It follows from Lemma \ref{lemWsujmart} and Doob's convergence theorem that for any vertex $u\in\etiq_{0}$ with generation at least $\underline{j}$ and any real $s\geq 0$,
\[
\as \qquad W_{s,u,j} \xrightarrow[j\to\infty]{} W_{s,u,\infty}\in [0,\infty).
\]

We can now establish the desired upper bound on the Hausdorff dimension of the sets $K_{u}$.

\begin{prp}\label{majdimKu}
For any vertex $u\in\etiq_{0}$ with $\gene{u}\geq\underline{j}$,
\[
\as \qquad \dim K_{u}\leq d_{*}.
\]
\end{prp}

\begin{proof}
Let us first assume that $d_{*}<0$. In particular, $\rho(0)$ is negative. The definition (\ref{defrhodim}) of the function $\rho$, together with the observation that $\alpha_{0,\gene{u}+\ell}=\ph_{1,\gene{u}+\ell}'(1)$ for any $\ell\geq 0$, implies that
\[
\prod_{\ell=0}^{j-1}\ph_{1,{\gene{u}}+\ell}'(1)\leq e^{\rho(0)j/2}
\]
for infinitely many integers $j\geq 1$. Hence, the number $\underline{\sigma}_{\gene{u}}$ defined by (\ref{defsigmaj}) is infinite and Lemma \ref{lemrecencfj} guarantees that $f_{\gene{u}}$ is equal to one. It follows from Lemmas \ref{lemlienfj}, \ref{equiKuvideZujext} and \ref{prplienGW} that the set $K_{u}$ is empty with probability one, so that its Hausdorff dimension is necessarily less than $d_{*}$.

Let us now assume that $d_{*}\geq 0$ and let us consider a real $s>d_{*}$. Then, $\rho(s)$ is negative, so that the limit inferior of $\prod_{\ell=\gene{u}}^{j-1} \alpha_{s,\ell}$ vanishes as $j\to\infty$. Meanwhile, $W_{s,u,j}$ converges almost surely to a finite limit. It follows that the limit inferior of $Z_{s,u,j}$ vanishes with probability one. Lemma \ref{tildeEsuleqliminfZsuj} ensures that $\tilde E_{s,u}$ vanishes almost surely. Lemma \ref{prplienEusHsKu} then implies that with probability one, $\dim K_{u}\leq s$. To deduce that the Hausdorff dimension is almost surely at most $d_{*}$, it suffices to let $s$ tend to $d_{*}$ along a decreasing sequence of reals.
\end{proof}

\section{Lower bound on the dimension}\label{arbresminor}

With a view to proving Theorem \ref{loidimTheta}, we establish in this section that the Hausdorff dimension of the sets $K_{u}$ defined by (\ref{defKu}) and composing $\Theta$ is almost surely at least the number $d_{*}$ defined by (\ref{defdetoile}), when they are nonempty. We may obviously assume that $d_{*}$ is positive. In particular, the number $\underline{j}$ given by (\ref{defunderj}) is finite. More precisely, we establish the following result.

\begin{prp}\label{prpdimKuW0uinfty}
Let us assume that $d_{*}$ is positive. Then, for any vertex $u\in\etiq_{0}$ with $\gene{u}\geq\underline{j}$, with probability one,
\[
K_{u}\neq\emptyset \qquad\Longrightarrow\qquad \dim K_{u}\geq d_{*}.
\]
\end{prp}

The rest of this section is devoted to the proof of Proposition \ref{prpdimKuW0uinfty}. Let us suppose that $d_{*}$ is positive and let us consider a vertex $u\in\etiq_{0}$ with generation at least $\underline{j}$. Our proof makes use of the processes $(\mathfrak{E}_{s,u,j})_{j\geq\gene{u}}$ defined by
\[
\forall s\in (0,d_{*}) \quad \forall j\geq\gene{u} \qquad \mathfrak{E}_{s,u,j} = {e_{s,j}}^{Z_{u,j}},
\]
where the reals $e_{s,j}$ are the probabilities defined by (\ref{defesj}).

\begin{lem}\label{lemfrakesujmart}
For any real $s\in (0,d_{*})$, the process $(\mathfrak{E}_{s,u,j})_{j\geq\gene{u}}$ is a nonnegative martingale.
\end{lem}

\begin{proof}
Let us consider an integer $j\geq\gene{u}$. Thanks to (\ref{Zuj1sum}), the conditional expectation of $\mathfrak{E}_{s,u,j+1}$, conditionally on $\trib_{j}$ is $\ph_{1,j}(e_{s,j+1})^{Z_{u,j}}$, which is equal to $\mathfrak{E}_{s,u,j}$ by virtue of Lemma \ref{lemrecesj}. To conclude, it suffices to observe that the variables $\mathfrak{E}_{s,u,\gene{u}},\ldots,\mathfrak{E}_{s,u,j}$ are $\trib_{j}$-measurable.
\end{proof}

It follows from Lemma \ref{lemfrakesujmart} and Doob's convergence theorem that
\[
\as \qquad \mathfrak{E}_{s,u,j} \xrightarrow[j\to\infty]{} \mathfrak{E}_{s,u,\infty}\in [0,1].
\]
The next lemma supplies a connection between the value of the limiting variable $\mathfrak{E}_{s,u,\infty}$ and that of the limit $W_{0,u,\infty}$ of the process $(W_{0,u,j})_{j\geq\gene{u}}$ defined by (\ref{defWsuj}).

\begin{lem}\label{W0uinftyposEsuinftynul}
For any real $s\in (0,d_{*})$,
\[
W_{0,u,\infty}>0 \qquad\Longrightarrow\qquad \mathfrak{E}_{s,u,\infty}=0.
\]
\end{lem}

\begin{proof}
Lemma \ref{lemmajesj} ensures that for any integer $j\geq\gene{u}$,
\begin{equation}\label{frakEsujleqZuj}
\mathfrak{E}_{s,u,j} \leq \left( 1 - \frac{1}{\varsigma_{s,j}} \right)^{Z_{u,j}} \leq \exp\left(-\frac{ Z_{u,j} }{ \varsigma_{s,j} }\right).
\end{equation}
Furthermore, note that $\phi_{s,\ell}'(1)=\alpha_{s,\ell}>0$ and $\phi_{s,\ell}''(1)\leq m_{\ell}\phi_{s,\ell}'(1)$ for any integer $\ell\geq\underline{j}$. As a result,
\[
\varsigma_{s,j}\leq\left( \prod_{\ell=\underline{j}}^{j-1} \phi_{s,\ell}'(1) \right)\sum_{n=j}^\infty\frac{ m_{n} }{ \prod\limits_{\ell=\underline{j}}^n \alpha_{s,\ell} }.
\]
Since $\rho(s)$ is positive and the sequence $(m_{j})_{j\geq 0}$ is bounded owing to Lemma \ref{lemnjborne}, we necessarily have, for $\eps\in (0,\rho(s))$ and $j$ large enough,
\[
\sum_{n=j}^\infty\frac{ m_{n} }{ \prod\limits_{\ell=\underline{j}}^n \alpha_{s,\ell} }\leq\sum_{n=j}^\infty e^{(\eps-\rho(s))(n+1)}=\frac{e^{(\eps-\rho(s))(j+1)}}{1-e^{\eps-\rho(s)}}.
\]
Letting $\eps=\rho(s)/2$, applying (\ref{defWsuj}) with $s=0$ so as to express $Z_{u,j}$ in terms of $W_{0,u,j}$ and observing that $0<\phi_{s,\ell}'(1)\leq\ph_{1,\ell}'(1)$ for any integer $\ell\geq\gene{u}$, we deduce that
\[
-\frac{Z_{u,j}}{\varsigma_{s,j}}\leq -\frac{ W_{0,u,j} }{\prod\limits_{\ell=\underline{j}}^{\gene{u}-1} \phi_{s,\ell}'(1)}(e^{\rho(s)/2}-1)e^{\rho(s)j/2}
\]
for all $j$ large enough. Therefore, if $W_{0,u,\infty}$ is positive, the right-hand side tends to $-\infty$ as $j\to\infty$, so that $\mathfrak{E}_{s,u,\infty}$ vanishes thanks to (\ref{frakEsujleqZuj}).
\end{proof}

To prove Proposition \ref{prpdimKuW0uinfty}, let us consider a real $s\in (0,d_{*})$. Due to Lemma \ref{lemfrakesujmart},
\[
\esp[\mathfrak{E}_{s,u,\infty}]=\esp[\mathfrak{E}_{s,u,\gene{u}}]=\prob(X_{u}=0)+e_{s,\gene{u}}\prob(X_{u}=1).
\]
Moreover, recall that if $X_{u}=0$, then the flow $E_{s,u}$ vanishes and that if $X_{u}=1$ with positive probability, then conditionally on the event $\{X_{u}=1\}$, this flow has the same distribution as the variable given by (\ref{Esuegalenloi}). Thus, the probability that $E_{s,u}$ vanishes is equal to the right-hand side of the previous equality. Therefore, this probability is equal to the expectation of $\mathfrak{E}_{s,u,\infty}$. In addition, Lemma \ref{prplienEusHsKu} shows that $E_{s,u}$ vanishes if the dimension of the set $K_{u}$ is less than $s$. As a consequence,
\[
\prob(\dim K_{u}<s)\leq\prob(E_{s,u}=0)=\esp[\mathfrak{E}_{s,u,\infty}].
\]
This last expectation may be written on the form
\[
\esp[\ind_{\{W_{0,u,\infty}>0\}}\ind_{\{\mathfrak{E}_{s,u,\infty}>0\}}\mathfrak{E}_{s,u,\infty}]+\esp[\ind_{\{W_{0,u,\infty}=0\}}\mathfrak{E}_{s,u,\infty}].
\]
The first term is at most the probability that $W_{0,u,\infty}$ and $\mathfrak{E}_{s,u,\infty}$ are both positive, which is equal to zero because of Lemma \ref{W0uinftyposEsuinftynul}, and the second term is at most the probability that $W_{0,u,\infty}$ vanishes. Furthermore, Lemmas \ref{lemlienfj} and \ref{prplienGW} imply that $W_{0,u,\infty}$ vanishes with probability
\begin{equation*}\begin{split}
&\,\prob(X_{u}=0)+\prob(\calW_{\gene{u},\infty}=0)\,\prob(X_{u}=1) \\
=&\,\prob(X_{u}=0)+\prob(\calZ_{\gene{u},j}\to 0 \text{ as } j\to\infty)\,\prob(X_{u}=1).
\end{split}\end{equation*}
Because of Lemmas \ref{lemlienfj}, \ref{equiKuvideZujext} and \ref{prplienGW}, this probability is also equal to that of the event $\{K_{u}=\emptyset\}$. We deduce that $\prob(\dim K_{u}<s)\leq\prob(K_{u}=\emptyset)$. Since the event $\{\dim K_{u}<d_{*}\}$ is the increasing union over $s\in (0,d_{*})$ of the events $\{\dim K_{u}<s\}$, this yields
\[
\prob(\dim K_{u}<d_{*})\leq\prob(K_{u}=\emptyset).
\]
Proposition \ref{prpdimKuW0uinfty} follows directly.

\section{Proofs of the main results}\label{arbresmainproofs}

\subsection{Proof of Theorem \ref{loidimTheta}}

Propositions \ref{majdimKueasy} and \ref{majdimKu} ensure that for any vertex $u\in\etiq_{0}$, the set $K_{u}$ defined by (\ref{defKu}) has Hausdorff dimension at most $d_{*}$ with probability one. Recall that, by virtue of its definition (\ref{defTheta}), the set $\Theta$ is the union over $u\in\etiq_{0}$ of the sets $K_{u}$. Hence, with probability one, the dimension of $\Theta$ is at most $d_{*}$. In particular, if $d_{*}$ is negative, then the dimension of $\Theta$ is almost surely equal to $-\infty$.

It remains to prove that if $d_{*}$ is nonnegative, then with probability one,
\[
\Theta\neq\emptyset \qquad\Longrightarrow\qquad \dim\Theta\geq d_{*}.
\]
We may obviously assume that $d_{*}$ is positive. If $\Theta$ is nonempty and has Hausdorff dimension less than $d_{*}$, then there exists a vertex $u\in\etiq_{0}$ such that $K_{u}\neq\emptyset$ and $\dim K_{u}<d_{*}$. By virtue of Proposition \ref{majdimKueasy}, the generation of $u$ is at least $\underline{j}$. Proposition \ref{prpdimKuW0uinfty} then ensures that such a vertex $u$ may exist only with probability zero. The result follows.

\subsection{Proof of Theorem \ref{thmprobThetavide}}

Lemma \ref{equiKuvideZujext}, along with (\ref{defTheta}), ensures that $\Theta$ is empty if and only if the processes $(Z_{u,j})_{j\geq\gene{u}}$, for $u\in\etiq_{0}$, become extinct. Moreover, for any integer $j\geq 0$, the processes $(Z_{u,n})_{n\geq\gene{u}}$, for $\gene{u}\leq j$, become extinct, if and only if the processes $(Z_{u,n})_{n\geq\gene{u}}$, for $u\in S_{j}$, do. Therefore,
\[
\prob(\Theta=\emptyset)=\lim_{j\uparrow\infty}\downarrow \prob(\forall u\in S_{j} \quad Z_{u,n}\to 0\text{ as }n\to\infty).
\]
In addition, for any $j\geq 0$, Lemma \ref{prplienGW} and the Markov condition (\ref{condMarkov}) imply that
\[
\prob(\forall u\in S_{j} \quad Z_{u,n}\to 0\text{ as }n\to\infty \:|\: X_{v}, \gene{v}\leq j)={f_{j}}^{\# S_{j}}.
\]
To conclude, it suffices to take expectations and to let $j\to\infty$.

\subsection{Proof of Proposition \ref{altthmprobThetavide}}

Let us suppose that $d_{*}$ is nonnegative and that $f_{j_{*}}$ vanishes for some integer $j_{*}\geq 0$. It follows from Lemma \ref{lemrecencfj} that $f_{j}=0$ and $\ph_{1,j}(0)=0$ for any integer $j\geq j_{*}$. Owing to (\ref{Phijrec}),
\[
\forall j\geq j_{*} \qquad \Phi_{j+1}(f_{j+1})=\Phi_{j}(f_{j})\cdot\ph_{0,j}(0)^{m_{0}\cdot\ldots\cdot m_{j-1}}.
\]
We conclude by arguing by induction and using Theorem \ref{thmprobThetavide}.

\subsection{Proof of Proposition \ref{condnecsufThetavide1}}

Let us assume that (\ref{conddimTheta1}) holds. For any integer $j\geq 0$, it follows from (\ref{recfctgeneSj}) that
\[
\esp[ {f_{j+1}}^{\# S_{j+1}} \:|\: X_{u},\ \gene{u}\leq j ] = \ph_{1,j}(f_{j+1})^{\# S_{j}}\ph_{0,j}(f_{j+1})^{m_{0}\cdot\ldots\cdot m_{j-1}-\# S_{j}}.
\]
In addition, Lemma \ref{lemrecencfj} ensures that $\ph_{1,j}(f_{j+1})=f_{j}$ and (\ref{condfjpos}) imply that $\ph_{0,j}(f_{j+1})$ belongs to the interval $(0,1]$. Taking expectations, we deduce that
\[
\Phi_{j+1}(f_{j+1}) \geq \Phi_{j}(f_{j})\cdot\ph_{0,j}(f_{j+1})^{m_{0}\cdot\ldots\cdot m_{j-1}}.
\]
Theorem \ref{thmprobThetavide} and Lemma \ref{lemrecencfj} then imply that for any integer $j_{0}\geq 0$,
\begin{equation}\label{minprThetavideprodj0}
\prob(\Theta=\emptyset) \geq \Phi_{j_{0}}(f_{j_{0}})\cdot\prod_{j=j_{0}}^\infty \ph_{0,j}\left(1-\frac{1}{\underline{\sigma}_{j+1}}\right)^{m_{0}\cdot\ldots\cdot m_{j-1}}.
\end{equation}
Note that $\Phi_{j_{0}}(f_{j_{0}})$ and the product above are both positive, owing to (\ref{condfjpos}) and (\ref{conddimTheta1}) respectively. Therefore, the set $\Theta$ is empty with positive probability.

Let us suppose that (\ref{conddimTheta1bis}) does not hold. For any integer $j\geq 0$, let $\overline{m}_{j}=m_{0}\cdot\ldots\cdot m_{j-1}$. Owing to (\ref{condfjpos}), the reals $f_{j}$ and $\ph_{0,j}(f_{j+1})$ belong to the interval $(0,1]$. Thus, by virtue of (\ref{recfctgeneSj}) and Lemma \ref{lemrecencfj},
\begin{equation*}\begin{split}
\esp[{f_{j+1}}^{\# S_{j+1}}\:|\: X_{u},\,\gene{u}\leq j] &= {f_{j}}^{\# S_{j}} \ph_{0,j}(f_{j+1})^{\overline{m}_{j}-\# S_{j}} \\
&\leq {f_{j}}^{\# S_{j}} \left( \ph_{0,j}(f_{j+1})^{\overline{m}_{j}/2}\ind_{\{\# S_{j}\leq \overline{m}_{j}/2\}} + \ind_{\{\# S_{j}> \overline{m}_{j}/2\}} \right) \\
&\leq {f_{j}}^{\# S_{j}}\ph_{0,j}(f_{j+1})^{\overline{m}_{j}/2}+{f_{j}}^{\overline{m}_{j}/2}.
\end{split}\end{equation*}
Taking expectations and then arguing by induction on $j$, one easily checks that
\begin{equation}\label{majespwjSj}
\forall j\geq 0 \qquad \Phi_{j}(f_{j})\leq\frac{f_{0}}{u_{j}}+\frac{1}{u_{j}}\sum_{k=1}^j u_{k}\,{f_{k-1}}^{\overline{m}_{k-1}/2},
\end{equation}
where, for any integer $j\geq 0$,
\[
u_{j}=\frac{1}{\prod\limits_{\ell=0}^{j-1}\ph_{0,\ell}(f_{\ell+1})^{\overline{m}_{\ell}/2}}\geq\frac{1}{\prod\limits_{\ell=0}^{j-1}\ph_{0,\ell}\left( 1-\frac{1}{\overline{\sigma}_{\ell+1}} \right)^{\overline{m}_{\ell}/2}}.
\]
Observe that $\sum_{j} {f_{j}}^{\overline{m}_{j}/2}<\infty$. Indeed, $\rho(0)$ is positive, $\ph_{1,\ell}''(1)\leq m_{\ell}\ph_{1,\ell}'(1)$ for any integer $\ell\geq\underline{j}$ and the sequence $(m_{j})_{j\geq 0}$ is bounded owing to Lemma \ref{lemnjborne}. As a result, we necessarily have, for $\eps\in (0,\rho(0))$ and $j$ large enough,
\[
\overline{\sigma}_{j}\leq\sum_{n=j}^\infty\frac{m_{n}}{\prod\limits_{\ell=j}^n\ph_{1,\ell}'(1)}=\left(\prod_{\ell=\underline{j}}^{j-1}\ph_{1,\ell}'(1)\right)\frac{e^{(\eps-\rho(0))(j+1)}}{1-e^{\eps-\rho(0)}}.
\]
Letting $\eps=\rho(0)/2$ and using Lemma \ref{lemrecencfj}, we obtain
\[
{f_{j}}^{\overline{m}_{j}/2}\leq\exp\left(-\frac{\overline{m}_{j}}{2\overline{\sigma}_{j}}\right)\leq\exp\left(-(e^{\rho(0)/2}-1)\frac{\overline{m}_{\underline{j}}}{2} e^{\rho(0)j/2}\right)
\]
for all $j$ large enough, which ensures the convergence of $\sum_{j} {f_{j}}^{\overline{m}_{j}/2}$. Moreover, the sequence $(u_{j})_{j\geq 0}$ is nondecreasing and diverges to infinity. Kronecker's lemma then guarantees that the right-hand side of (\ref{majespwjSj}) tends to zero as $j\to\infty$, see \cite[p. 103]{Dienes:1957fk}. Theorem \ref{thmprobThetavide} implies that $\Theta$ is empty with probability zero.

\subsection{Proof of Proposition \ref{condnecsufThetavide2}}

Let us assume that (\ref{conddimTheta2}) holds. If $\underline{\sigma}_{\underline{j}}$ is infinite, then one easily checks using (\ref{defsigmaj}) that $\underline{\sigma}_{j}=\infty$ for any $j\geq\underline{j}$. In particular, $f_{\underline{j}}=1$ thanks to Lemma \ref{lemrecencfj} and $\ph_{0,j}(1-1/\underline{\sigma}_{j+1})=1$ for each $j\geq\underline{j}$. Applying (\ref{minprThetavideprodj0}) with $j_{0}=\underline{j}$, we deduce that $\Theta$ is empty with probability one. Conversely, if $\underline{\sigma}_{\underline{j}}$ is finite, then the function $\Phi_{\underline{j}}$ and the functions $\ph_{0,j}$, for $j\geq\underline{j}$, are constant and equal to one. Together with (\ref{minprThetavideprodj0}), this directly implies that $\Theta$ is almost surely empty.

Let us suppose that $\Theta$ is empty with probability one. By Theorem \ref{thmprobThetavide}, the expectation of ${f_{j}}^{\# S_{j}}$ is equal to one for any integer $j\geq 0$. Let us assume that $\overline{\sigma}_{\underline{j}}$ is finite. Using (\ref{defsigmaj}), one easily checks that $\overline{\sigma}_{j}$ is also finite, for any $j\geq\underline{j}$. By Lemma \ref{lemrecencfj}, the probability $f_{j}$ is less than one for each $j\geq\underline{j}$, so that with probability one, $\# S_{j}$ vanishes for all $j\geq\underline{j}$. In particular, $\Phi_{\underline{j}}(0)=1$ and $\ph_{0,j}(0)=1$ for any $j\geq\underline{j}$, thanks to (\ref{recfctgeneSj}). Thus, (\ref{conddimTheta3}) necessarily holds.

\section{A straightforward generalization}\label{arbresgeneralisation}

In this last section, the families $(L_{u})_{u\in\etiq_{0}^*}$ and $(X_{u})_{u\in\etiq_{0}}$ are not assumed to be independent anymore and Assertions (\ref{condcompact3}) and (\ref{condMarkov}) are replaced by the following one:
\renewcommand{\theenumi}{\Alph{enumi}}
\begin{enumerate}\setcounter{enumi}{6}
\item For any $u\in\etiq_{0}$, the conditional law of $(X_{u1},\ldots,X_{um_{\gene{u}}},L_{u1},\ldots,L_{um_{\gene{u}}})$, conditionally on the variables $X_{v}$ and $L_{v}$, for $v\not\in u\etiq_{\gene{u}}^*$, is $\lambda_{X_{u},\gene{u}}$.
\end{enumerate}
Here, $\lambda_{t,j}$ denotes a fixed probability measure on $\{0,1\}^{m_{j}}\times [\underline{\beta},\overline{\beta}]^{m_{j}}$ for every $t\in\{0,1\}$ and every $j\geq 0$. Note that the initial case in which Assertions (\ref{condcompact3}) and (\ref{condMarkov}) hold and the families $(L_{u})_{u\in\etiq_{0}^*}$ and $(X_{u})_{u\in\etiq_{0}}$ are independent may be recovered by letting $\lambda_{t,j}$ be the product measure $\nu_{t,j}\otimes\mu_{j}$, for any $t$ and any $j$.

Then, it is relatively straightforward to adapt the proofs exposed in Sections \ref{arbresprel}-\ref{arbresmainproofs} above in order to establish that Theorems \ref{loidimTheta} and \ref{thmprobThetavide}, as well as Propositions \ref{altthmprobThetavide}, \ref{condnecsufThetavide1} and \ref{condnecsufThetavide2}, still hold this generalized context, provided that (\ref{defalphasj}) and (\ref{defphtj}) are replaced by
\begin{equation*}\begin{split}
\alpha_{s,j} &=\int_{\{0,1\}^{m_{j}}\times [\underline{\beta},\overline{\beta}]^{m_{j}}} \sum_{k=1}^{m_{j}} {\ell_{k}}^s x_{k} \ \lambda_{1,j}(\dd x,\dd\ell) \\
\text{and}\qquad \ph_{t,j}(z) &=\int_{\{0,1\}^{m_{j}}\times [\underline{\beta},\overline{\beta}]^{m_{j}}} z^{x_{1}+\ldots+x_{m_{j}}} \lambda_{t,j}(\dd x,\dd\ell),
\end{split}\end{equation*}
respectively.

\end{document}